\documentclass{ifcolog-c}
\usepackage{amsmath,amsfonts}
\usepackage[utf8]{inputenc}
\usepackage[german,english]{babel}
\usepackage[authoryear]{natbib}
\let\parencite\citep
\usepackage{doi}
\hypersetup{colorlinks,allcolors=blue}

\newenvironment{axiom}[1]
    {\begin{list}{}%
    {%
    \settowidth\labelwidth{#1}%
    \setlength\leftmargin{\labelwidth+\labelsep}}}%
    {\end{list}}

\title{Some Unpublished Letters by Gödel and von Neumann in the Fraenkel Archive}

\titlethanks{We gratefully acknowledge permission to publish the
letters from Gödel and von Neumann to Fraenkel by the National Library
of Israel, the Institute for Advanced Study, Princeton (for the
letters by Gödel), and Marina von Neumann Whitman (for the letters by
von Neumann). Thanks to Shaul Greenstein, archivist at the National
Library of Israel, for his assistance in locating the correspondence
between Fraenkel, Gödel, and von Neumann. Thanks also to Juliette
Kennedy, Mark van Atten, Paola Cantù, Ulf Hashagen, and especially Aki
Kanamori.}
\addauthor[mancosu@socrates.berkeley.edu]{Paolo Mancosu}{University of California, Berkeley}
\authorrunning{Mancosu}
\addauthor[rzach@ucalgary.ca]{Richard Zach}{University of Calgary}
\authorrunning{Mancosu \& Zach}

\titlerunning{Some Unpublished Letters in the Fraenkel Archive}

\def\hbreak{\centerline{\rule{.5\textwidth}{.5pt}}}
\def\sout#1{$\langle\!\langle$#1$\rangle\!\rangle$}
\newcommand{\inter}[1]{$[\![$#1$]\!]$}
\def\pag#1{$\mid^{#1}$}
\newcommand{\closing}[1]{\par\smallskip\begin{raggedright}\parindent=0pt\setlength{\parskip}{.5\baselineskip}
#1\end{raggedright}}
\newcommand{\opening}[1]{\par\begin{flushleft}#1\end{flushleft}\smallskip\par\noindent\ignorespaces}
\newcommand{\return}[1]{\begin{flushright}#1\end{flushright}}

\begin{document}

\selectlanguage{english}

\maketitle

\section{Introduction}

In the present paper we edit and discuss several letters by Kurt Gödel
and Johann (János) von Neumann which are found in the archive of
Abraham Adolf Fraenkel at the National Library of Israel. The archive
is still uncatalogued.

Adolf (later Abraham) Fraenkel (1891--1975) was a professor of
mathematics at Marburg University until 1928. In 1928 he moved to Kiel
and after one year he accepted a position at the Hebrew University of
Jerusalem where he remained for the rest of his professional career
(for a biographical memoir see \citealt{Fraenkel2016}). The existence
of some of the letters we publish was mentioned by Ivor
Grattan-Guinness, who had been shown them by Fraenkel's widow when the
letters where still in her possession. Referring to ``an astounding
14-page letter dated 26 October 1923'' by von Neumann, published
below, Grattan-Guinness wrote:
\begin{quote}
    I saw this letter in the presence of his widow in 1982, at her
    flat in Jerusalem; it is now in his Papers. It includes also
    letters from Carnap, Gödel and Quine. No copy of this letter
    exists in the von Neumann Papers, held at the Library of Congress,
    Washington. \citep[note 33, p. 479]{Grattan-Guinness2000}
\end{quote}
As mentioned, we publish here all the letters from Gödel and von
Neumann contained in the Fraenkel archive. While some are more
valuable than others, we include everything we found in the spirit of
contributing to making available more unpublished correspondence by
Gödel and von Neumann and, consequently, to contribute to a better
understanding of some crucial moments in the development of the
foundations of mathematics in the 1920s and 1930s.

Both Gödel and von Neumann are towering intellectual figures of the
20th century. Their importance is highlighted by two trade biographies
published recently and making the bestseller lists
\citep{Bhattacharya2021,Budiansky2021}. Academic interest in Gödel's
thought also continues to grow, fueled by the outstanding achievement
of the five volumes of Gödel's papers and correspondence (see
\citealt{Godel2003,Godel2003a} for the correspondence). More recent
archival work continues to shed light on the development of Gödel's
life and thought.\footnote{See, among other publications for the last
ten years, \citet{AdzicDosen2016}, \citet{CroccoEngelen2016} and
several articles therein contained, \citet{Godel2019},
\citet{KanckosLethen2021}, \citet{Kennedy2020}, \citet{Lethen2021},
\citet{vanAtten2015}, and \citet{vonPlato2020}. Transcription from the
Philosophical notebooks by Gabriella Crocco's research group are
published on HAL at the following link:
\url{https://hal.archives-ouvertes.fr/search/index/q/*/authIdHal_s/gabriella-crocco}}
Von Neumann's early biography is also receiving its well-deserved
attention by historians; see in particular
\citet{Hashagen2006,Hashagen2006a,Hashagen2010}.

We think it fitting to present this paper in a volume honoring the
scholarly achievements of Volker Peckhaus, whose archival work has
dramatically increased our understanding of the history of the
foundations of logic and mathematics and whose personal assistance
during our own investigations on these subjects has been immensely
valuable to us. As we will see, the correspondence with von Neumann
concerns a pivotal moment in the history of set theory, the study of
which is one of Volker's main interests \citep{Peckhaus1990a,Peckhaus1992,Peckhaus2005,EbbinghausPeckhaus2015}.

\section{Letters from Gödel}

The letters from Gödel found in the Fraenkel archive fall into two
distinct groups.\footnote{No letters between Fraenkel and Gödel are
listed in the finding aid of the Gödel Papers published in
\citet[pp.~469--562]{Godel2003a}.} To the first group belong three
fragments dating from 1931--1932. The longer fragment (four
handwritten pages) comes from a letter datable to the end of 1931 and
corresponds, as can be gathered from the numbering on the upper left
corner of the extant pages, to the last four pages of a letter whose
first two pages are missing. The other two, much smaller, fragments
come from the last two pages of a letter datable to 1932. In this case
it is not possible to know how long the original letter was. Only
the lower halves of the last two pages, from which the fragments
originate, are extant. These fragments are historically the more
interesting, since they concern the relationship between Gödel's
incompleteness theorem \citep{Godel1931} and work by Herbrand,
Presburger, and Zermelo.

To the second group belong three complete letters written in 1958. We
will first briefly describe this second group of letters, which is
merely of biographical interest, and then devote more attention to the
fragments from 1931--1932. The first letter in this second group
concerns Michael Rabin and Azriel Levy. Both Rabin and Levy had
studied with Fraenkel at the University of Jerusalem. Rabin had done
his Ph.D. at Princeton in 1956 under Alonzo Church; Levy obtained his
Ph.D. in 1958 at the University of Jerusalem under Fraenkel and
Abraham Robinson. Fraenkel had apparently asked Gödel's opinion in
regards to a possible appointment for Rabin as Assistant Professor in
one of the previous letters Gödel refers to, of 2~December 1957 and
2~January 1958. Gödel provides a very positive evaluation of Rabin's
work, but says he can't comment on his abilities as a teacher. In the
second part, Gödel responds to an inquiry by Fraenkel concerning a
possible fellowship for Levy at the Institute for Advanced Study. He
writes that it is too late in the year to make arrangements. In the
letter of 12~May 1958, Gödel mentions funds have become available, but
that he has learned that Levy has found another position in the
meantime, and so would not be able to take up the fellowship at the
IAS. He suggests the possibility of renewing the offer of a fellowship
for the 1959/60 academic year. The final letter concerns a visit by
Fraenkel to Princeton. Gödel mentions that he usually does not attend
the department tea, and asks Fraenkel to make arrangements by
telephone.

We now move to the first group of (fragments of) letters from
1931--1932, which are of greater scientific interest. The first
fragment begins with Gödel commenting on someone's work in set theory,
as if he had been asked something about it by Fraenkel in a
(presumably lost) previous letter. We do not know which paper or book
Gödel may be referring to here. The second, long paragraph concerns
\cite{Herbrand1931}. Fraenkel reviewed \cite{Herbrand1931} for the
\emph{Jahrbuch für die Fortschritte der Mathematik} (JFM 57.0056.04).
We have no other letter by Gödel where he discusses Herbrand's article
on the consistency of arithmetic at length. However, there are brief
mentions: in a letter to Heyting of 15 November 1932
\citep[p.~61]{Godel2003a} and in an exchange with van Heijenoort (14
August, 27 August, 18 September 1964) \citep[pp.~316,
318--319]{Godel2003a}. These mentions occur much later and in the
context of interpreting what Herbrand meant by strict intuitionism. It
is true that Gödel's letter to Herbrand dated 25~July 1931
\cite[pp.~20--21]{Godel2003a} addresses some of the issues under
consideration. The reason is that Herbrand wrote to Gödel on 7~April
1931, in order to discuss how Gödel's results fit with results he had
obtained previously (the nomenclature of the theories he discusses are
the same as those of \citealt{Herbrand1931}).

Herbrand's article \citeyearpar{Herbrand1931} was completed on 14~July
1931. It was received by the journal on 27~July 1931 and appeared in
the first issue of volume 166 (1932), which was released on
12~November 1931. Since the second \emph{Heft} came out on 10~December
1931, and Gödel refers to the ``latest'' issue of \emph{Journal für
die Reine und Angewandte Mathematik} (which was usually referred to by
its subtitle, \emph{Crelle's Journal}), we can date this letter
approximately as written between mid-November and mid-December 1931,
allowing for a possible shift forward depending on how long it would
have taken for the issue to reach Gödel or the library used by Gödel.
There is however something puzzling in the letter. This is the final
request concerning von Neumann: ``Should you, dear professor, have
corresponded with Mr.\ v.~Neumann about the question of the
consistency proof, I would be grateful if you could communicate to me
his views.'' What is surprising here is the fact that Gödel knew
perfectly well what von Neumann's position was since they had
corresponded on the issue and he had mentioned what von Neumann's
position was already in the presentation of his results in front of
the Schlick Circle in January 1931 (see \citealt{Mancosu1999a}). Von
Neumann had expressed his position on the incompleteness theorems in
the letter to Gödel dated 10~January 1931
\citep[pp.~338--339]{Godel2003a}. But perhaps Gödel wanted to know
whether von Neumann's position had shifted.

In his letter to Fraenkel, Gödel indicates his acceptance of all of
the mathematical theorems demonstrated by Herbrand. In particular, he
agrees with Herbrand that there is no inconsistency between Herbrand's
consistency proof and Gödel's incompleteness theorems. Indeed, the
proof of consistency given by Herbrand, for a system weaker than
Weyl's system given in \emph{Das Kontinuum} \citep{Weyl1918}, is not
formalizable in the system itself even though the consistency proof is
unobjectionable from the intuitionistic, i.e., finitary,\footnote{At
the time, intuitionism and finitism were identified.} point of view.

Gödel's position concerning Herbrand's methodological standpoint
corresponds to the opinions he expressed also on other occasions (see
the correspondence with Bernays, Herbrand, and von Neumann published
in \citealt{Godel2003,Godel2003a}). Gödel does not agree with all of
Herbrand's conjectures. In particular, he discusses the issue of the
extent of “finitism”. He points out that Hilbert and his school use
the notion in a rather intuitive and contentful (i.e., meaningful) way
and that their finitist proofs are not formalized in a specific
system. Gödel thinks, in agreement with Herbrand and the
intuitionists, that a definition of “finitism” cannot be given in
principle. Where he disagrees with Herbrand and von Neumann is in the
claim that all finitary proofs can be formalized within the axiom
system of classical analysis. Herbrand had made the following
conjecture in a letter to Gödel (and in almost identical words in
\citealt{Herbrand1931}):
\begin{quote}
    It reinforces my conviction that it is impossible to prove that
    every intuitionistic proof is formalizable in Russell's system,
    but that a counterexample will never be found. There we shall
    perhaps be compelled to adopt a kind of logical postulate.
    \citep[p.~21]{Godel2003}
\end{quote}
Gödel, by contrast, believed that it could be in principle possible
that some finitary proofs could go beyond the formal proof
procedures of classical analysis and thus, contrary to Herbrand and
von Neumann, he did not believe that his result showed the
impossibility of carrying out Hilbert's program. This was the position
he had already put forth in his 1931 article. Gödel in his 1931 paper
states that his results 
\begin{quote}
    do not contradict Hilbert's formalistic viewpoint. For this
    viewpoint presupposes only the existence of a consistency proof in
    which nothing but finitary means of proof are used, and it is
    conceivable that there exist finitary proofs that cannot be
    expressed in the formalism of~$P$.
    \cite[p.~195]{Godel1986}\footnote{$P$ refers to the modification
    of the system of \emph{Principia Mathematica} used, and proved
    incomplete, in \citet{Godel1931}.}
\end{quote}

In the letter to Fraenkel, Gödel grants that the finitary arguments so
far used in the Hilbert's school do not suffice to prove the
consistency of analysis (or arithmetic). Gödel, however, continued to
take the view taken in his 1931 article, that the possibility remains
open that finitary modes of reasoning capable of establishing this
consistency result might still be found and employed to carry out
Hilbert's program. Gödel also objects to Herbrand's conjecture by
claiming that it is almost meaningless or unintelligible, for it would
require a definition of “finitary proof” for it to make sense (a
similar claim is made in his letter to Herbrand dated 7~April
1931).\footnote{“I would like now to enter into the question of the
formalizability of intuitionistic proofs in certain formal systems
(say that of Principia mathematica), since here there appears to be a
difference of opinion. I think, insofar as this question admits a
precise meaning at all (due to the undefinability of the notion
“finitary proof”, that could justly be doubted), the only correct
standpoint can be that we admit not knowing anything about it”
\cite[p. 23]{Godel2003a}.}

As pointed out, Gödel's position on finitism in this letter
corresponds to that defended in other letters from the period. Gödel
will later hold, starting in 1933 (see \citealt[p.~8]{Sieg2003}), that his
incompleteness theorems show that Hilbert's program, as originally
conceived, cannot be carried out. For a more detailed discussion of
Gödel's evolution with respect to the notion of finitism we refer to
\cite{Feferman2003}, \citet{Sieg2003}, and \citet{Mancosu2004}.

Let us conclude with the remaining two fragments which, as mentioned
earlier, stem from a single letter. In the first fragment of this
second letter, Gödel stresses that the article by
\cite{Presburger1930} does not conflict with his own incompleteness
result. This is the only place in Gödel's entire correspondence where
Gödel cites Presburger's result. The contribution by Presburger was
cited in a letter from Bernays to Gödel dated 18 January 1931
\citep[p.~90]{Godel2003}. Gödel also refers to Presburger in his 1934
Princeton lectures (see \citealt[p.~367]{Godel1986}). Fraenkel also
reviewed \cite{Presburger1930}, in JFM 56.0825.04.

The second fragment of this second letter reads (almost) in full:
\begin{quote}
    I would like to mention on this occasion that Zermelo in
    \emph{Jahresbericht} Bd~41, 2.~Abt. p.~87
    \inter{\citealt{Zermelo1932}} has made completely inadequate
    claims about my paper \inter{\citealt{Godel1931}}, to which I hope
    to respond in the \emph{Jahresbericht} in due course.
\end{quote}
The reference to Zermelo's article in the \emph{Jahresbericht} and
Gödel's indicated intention to reply, which he did not pursue, allows
us to safely date this letter to 1932.  For more on the Gödel-Zermelo
exchange see \citet{Grattan-Guinness1979}, and
\citet{Dawson1985,Dawson2003}.

\section{Letters from von Neumann}

The Fraenkel papers contain six letters from von Neumann.  The first
is the ``astounding'' 14-page letter mentioned by Grattan-Guinness,
dated 26~October 1923.  At the time, 19-year-old von Neumann was
nominally studying chemical engineering at the ETH
Zürich.\footnote{See \citet{Hashagen2006} on von Neumann's student
years.} He had just completed  two years of studying chemistry at the
University of Berlin (in preparation for the ETH entrance exam), as
well as mathematics and physics.  There, he had caught the eye of
Erhard Schmidt, who shared von Neumann's interest in set theory.  In
August 1923, von Neumann completed a manuscript on set theory, which
he sent to Zermelo (a lifelong friend of Schmidt's), and which Schmidt
sent to Fraenkel. Fraenkel later recalled:
\begin{quote}
    In 1923, Erhard Schmidt, who had been a professor in Berlin since
    1920 and was one of the editors of the journal \emph{Mathematische
    Zeitschrift}, sent me a lengthy manuscript while I was in Marburg
    entitled ``Die Axiomatisierung der Mengenlehre'' (The
    Axiomatization of Set Theory), written by Johann von Neumann, a
    name as yet unfamiliar to me. Schmidt, who a year earlier had
    presented a paper of mine to the Prussian Academy of Sciences,
    accompanied the manuscript with a letter saying that they could
    find no one in Berlin who was competent to evaluate the
    manuscript, and so he requested a statement from me whether a
    publication in the \emph{Mathematische Zeitschrift} should be
    considered. The article, 84 printed pages, did in fact appear five
    years later in the 27th volume of the journal, but had already
    been published in essentially the same form in 1925 in Hungarian
    as von Neumann's doctoral dissertation in Budapest. At around the
    same time, with an accompanying letter of August 14, 1923, von
    Neumann sent me his work directly, and expanded it two months
    later with comments in a letter comprising 14 large quarto-sized
    pages.
    \citep[143]{Fraenkel2016}
\end{quote}
The letter of 14 August appears to be lost, but it likely does not
differ much from the letter von Neumann sent to Zermelo (see
\citealt[pp.~271--273]{Meschkowski1967}). The first letter in the
Fraenkel archive is evidently a response to a reply by Fraenkel.

As Fraenkel acknowledges following the passage quoted above, ``it took
considerable effort for me to work through the treatise, which not
only deviated from everything that had thus far been published on the
axiomatization of set theory, but it also contradicted previous
ideas.''  The fundamental difference between von Neumann's approach
and those before him (such as Zermelo, Skolem, and Fraenkel himself)
concerns both the choice of primitives and the choice of axioms.
Rather than taking sets and classes as primitive, von Neumann's theory
distinguished between arguments (``\emph{I~Dinge}'', I~objects) and
functions (``\emph{II Dinge},'' II~objects). Objects which are both
functions and arguments are ``\emph{I~II Dinge}.'' The primitive
operations are $[x,y]$, the application of the function~$x$ to
argument~$y$, and $\langle x,y\rangle$, the ordered pair of $x$
and~$y$.\footnote{\Citet{vonNeumann1925} uses $(x,y)$ for ordered
pairs. For von Neumann's notations and definitions, see
\citet{vonNeumann1928}.} In this framework, sets (\emph{Mengen}) and
classes (\emph{Bereiche}) are represented by functions
(II~objects)~$a$ which take only one of two possible values, $A$
and~$B$, which are fixed I~objects: the corresponding class
contains~$x$ iff $[a, x] = B$, or equivalently: $[a, x] \neq A$. Thus
sets and classes, in von Neumann's system, are both II~objects. The
difference lies in whether they can themselves be elements of other
classes (sets can be elements; classes cannot be). So in von Neumann's
framework, the difference comes down to whether they can be
\emph{arguments}, i.e., I~objects. Something that is both a function
and an argument is a I~II object; so a set is a I~II object that takes
only $A$ and~$B$ as values, and a proper class is a II~object like
this which isn't a I~object.

The letter to Fraenkel sets out von Neumann's motivations for this
approach. The first reason is that, according to von Neumann, an
axiomatic set theory stands and falls with its ability to develop a
theory of ordinals. But Zermelo's axioms are not enough to do this.
Cantor's approach was to define ordinals as order types. This approach
is not available to Zermelo, for, as von Neumann puts it,
\begin{quote}
nothing guarantees \inter{in Zermelo's set theory}, that the ``ordinal
number'' of a set that is not ``too large'' doesn't itself become
``too large.'' After all, the only thing one knows about this ordinal
number is that it is similar to a given set (and thus equivalent).
From this one can only conclude that the ordinal number is not too
large if one is in possession of its ``replacement axiom.''
\end{quote}
The replacement axiom had been introduced independently by
\citet{Fraenkel1922a} and \cite{Skolem1923}.\footnote{For an in-depth
discussion of the replacement axiom and its history, see
\citet{Kanamori2012}.} Briefly, it states that if $X$ is a set and $f$
is a function mapping the elements of~$X$ to sets, then ``replacing''
every element~$a$ of $X$ by its value~$f(a)$ results in a set.
However, von Neumann continues,
\begin{quote}
In Zermelo's axiomatics it is hardly possible to accommodate the
replacement axiom. For the definition of ``function'' in Zermelo (set
of all pairs $x, f(x)$) requires from the outset that the domain as
well as the range are ``not too large.'' The replacement axiom, which
stipulates that the former implies the latter, would then indeed be
superfluous, and not even expressible in its strongest form
\inter{\emph{in voller Schärfe}}. One is therefore forced to introduce
the notion of function with arbitrary domain in a different way, or
independently.  And I think the absolutely simplest way to do this is
to assume the ``function'' with arbitrary domain \emph{and} range from
the outset.
\end{quote}
Von Neumann's system, however, does not include the replacement axiom.
Instead, it contains the central axiom IV\,2:
\begin{axiom}{IV\,2}
\item[IV\,2] A II~object~$a$ is not a I~II object if and only if
    there is a II~object~$b$ such that for every I~object~$x$ there is
    a~$y$ with $[a,y] \neq A$, $[b,y]=x$.
\end{axiom}
Recall that I~II objects correspond to sets, and II~objects that
aren't I~II objects to proper classes. So this axiom roughly says that
$a$ is a proper class iff there is a function~$b$ such that, for
any~$x$ whatsoever, there is  some ``element''~$y$ of~$a$ such that
the value of~$b$ for argument~$y$ is~$x$. In other words, proper
classes are exactly those classes for which there is a surjective map
to the universe.

This does not at all look like the replacement axiom, and from von
Neumann's reply we can assume that Fraenkel had asked for a
justification of axiom IV\,2. Von Neumann provides it on pp.~3--5.
First of all, he grants that one could replace it by axioms which
more directly correspond to replacement and choice, namely:
\begin{axiom}{III\,2$^*$}
    \item[IV\,2$^*$] Let $a, c$ be II objects, $b$ a I~II object.
    Suppose for every $y$ such that $[a,y] \neq A$ there is an $x$
    such that $[b,x] \neq A$, $[c, x] = y$. Then $a$ is also a I~II
    object.
    \item[III\,2$^*$] Let $a$ be a II object. Then  there is a II
    object~$b$, so that whenever there is a $y$ with $[a, \langle x,
    y\rangle] \neq A$, then $[b, x]$ is such a~$y$. 
\end{axiom}
Similar considerations can be found in the proof in
\citet[p.~230ff]{vonNeumann1929a} that the system of
\citet{vonNeumann1928} is equiconsistent with the system~$S^*$, which
replaces axiom IV\,2 by the above two.

However, what counts in favor of IV\,2 are (1) that it is a clear
definition of when a set is ``too large,'' (2)~it implies the axiom of
choice in the strong form of having a global choice function, and
(3)~it allows the development of a theory of finite sets without the
``infinite'' axioms V\,2 and~V\,3.

On p.~5, von Neumann compares axiom IV\,2 with the separation
axiom, which can be formulated in von Neumann's system as:
\begin{axiom}{IV\,2$^{**}$}
    \item[IV\,2$^{**}$] Let $a$ be a II object, and $b$ a I~II
    object. Suppose $[a, x] \neq A$ implies $[b, x] \neq A$. Then $a$
    is also a I~II object.
\end{axiom}
This follows from axiom~IV\,2 (cf. \citealt[p.~685]{vonNeumann1928}).
Conversely, axiom IV\,2 follows from IV\,2$^{**}$ together with the
requirement that the class~$\Omega$ of all I~objects has a
well-ordering~$W$ in which all proper initial segments of~$\Omega$ are
sets.\footnote{See \citet{Kanamori2009} on the relationship between
von Neumann's original axiomatization and that of Bernays, what is now
known as the von Neumann-Bernays-Gödel (NBG) system.}

Von Neumann goes on to discuss some further aspects of parts I and~II
of the manuscript, namely: 
\begin{enumerate}
\item a critique of axiom group V (about
infinity), 
\item the reasons for why equivalence (sameness of cardinality)
is treated after similarity (sameness of order type), 
\item the failure
of the claim, made in the manuscript, that not only well-ordered sets
but also well-ordered classes are pairwise comparable, 
\item a proof sketch (missing in the manuscript) of the possibility of
definition by transfinite induction (cf. \S IX of
\citealt{vonNeumann1928}).
\end{enumerate}

On pp.~9--11 of the letter, von Neumann discusses part~III of the
manuscript, which apparently concerned methodological reflections on
the purpose of an axiomatic development of set theory and a comparison
with an ``intuitionistic'' development of the same, as well as the
question of ``relativity'' of set-theoretic notions such as
(un)countability in the sense of \citet{Skolem1923}.\footnote{See
\citet{Jane2001} for discussion of Skolem's arguments.} Von Neumann
first clarifies that his approach takes the logical notions as
meaningfully given, whereas the notions of set, function, and choice
are treated axiomatically---these are the notions which are suspected
of being ``impredicative'' and so objectionable from the standpoint of
Brouwer and Weyl (cf. \S1 of \citealt{vonNeumann1925}). He states that
set theory can be integrated into Hilbert's axiomatization of
mathematics \citep{Hilbert1923}. He also indicates agreement with
\citet{Skolem1923} insofar as the ``countability of the system'' does
not conflict with the ``axiomatic existence of uncountable sets.'' (He
notes that he was not aware of Skolem's paper when composing the
manuscript and thanks Fraenkel for drawing his attention to it.) He
writes:
\begin{quote}
As regards the ``relativity,'' or as I put it, the ``irreality'' of
cardinalities, I am in agreement with Skolem. However, I do not draw
the same conclusions from it that he does: I believe, rather, that it
provides a strong argument \emph{against} intuitionism. One will never
even come close to the concept of an uncountable continuum on an
intuitionistic basis, since the actual world of ideas is countable.
The continuum is essentially transfinite, impredicative, and can
\emph{only} be captured by the formalist method.
\end{quote}
He continues to explain that, following Hilbert, such a formalist
treatment of the uncountable requires a proof of consistency of the
system. A consistency proof of the system containing logical
implication, negation, identity, and choice is easily
possible.\footnote{Here von Neumann presumably has in mind Hilbert's
consistency proof of logic in the $\varepsilon$-calculus in
\citet{Hilbert1923}. See \citet{Zach2003,Zach2004a} on Hilbert's
consistency proofs and the relation to \citet{vonNeumann1927}, and
\citet{Bellotti2016} on von Neumann's consistency proof specifically.}
To deal with set theory, the notion of function must also be
introduced, but for this system no consistency proof is as of yet
available. He would attempt to provide such a consistency proof in
\citet{vonNeumann1927}.

Finally (pp.~11--13), von Neumann discusses the question of
categoricity of the system and the relativity of cardinality and
well-ordering. Here, he raises the same points as he does in more
detail in \S II of \citet{vonNeumann1925}: the axiom system does not
settle whether there are I~objects that are not I~II objects and it
does not settle whether non-well-founded sets exist or not (e.g., sets
that satisfy $a = \{a\}$).\footnote{Von Neumann identifies these as
Mirimanoff's ``ensembles extraordinaires,'' which establishes that he
was familiar with \citet{Mirimanoff1917} already in 1923.} These
questions could be settled by additional axioms, but more importantly,
the relativity of well-orderedness and finitude suggests that set
theory is inherently non-categorical. Not only is there no maximal
model, but there need not be a minimal model either. Before he gives a
proof sketch of the relativity of well-orderedness, he writes:
\begin{quote}
    Following all this, I believe that set theory is incapable of
    being axiomatized categorically. Then, however, there is no
    categorical axiomatic system at all in the world (with the
    exception of finite systems). (Systems which seem to be
    categorical, such as mathematics, depend on set theory.) I
    am as of yet unsure of the scope and importance of this fact.
\end{quote}

Von Neumann concludes the letter by mentioning that Schmidt had
promised to publish the manuscript (after revisions which von Neumann
considered necessary), in the \emph{Mathematische Zeitschrift}. He
mentions that the article would run to 40 printed pages, and asks
Fraenkel if an earlier publication without these changes would be
advisable, and if so, in what venue. Finally, von Neumann wonders if
the theory of ordinals that he proposed is new.

In addition to the inherent interest of von Neumann's extensive
discussion, which gives us a glimpse of the early development of von
Neumann's theory, the letter also sheds light on the relationship
between von Neumann's publications on set theory.  It has been assumed
that the manuscript is more or less identical to
\citet{vonNeumann1928}. However, that paper runs to 84 printed pages,
over twice as long as the manuscript sent to Zermelo and Fraenkel. At
least the proof of transfinite induction therein was not yet included
in the manuscript.  We also see that many topics dealt with in
\citet{vonNeumann1925} (intuitive vs. axiomatic development of set
theory, models of set theory, Skolem's paradox, categoricity),
\citet{vonNeumann1927} (consistency proof of the notion of function),
and \citet{vonNeumann1929a} (relation between axiom IV\,2 on the one
hand, and replacement, separation, and choice on the other) were
already on von Neumann's mind in 1923.

The second letter is dated 17 June 1927.  At this point, von Neumann
was just completing his semester-long visit to G\"ottingen after
graduating in chemical engineering at the ETH in October 1926. He had
also just submitted his application for \emph{Habilitation} in
mathematics at the University of Berlin on 26 April 1927.\footnote{The
\emph{Habilitationsschrift} was entitled ``Der axiomatische Aufbau der
Mengenlehre,'' very likely a version of \cite{vonNeumann1928}. There
apparently was concern that the contents of this work (the version
submitted is not preserved), was too close to von Neumann's 1925
dissertation at the University of Budapest, and he ultimately replaced
it with a paper on operator theory. See
\citet{Hashagen2006,Hashagen2010}.} Von Neumann provides a proof of
the existence of the set~$\mathsf{M}$ ($\{\omega, \wp(\omega),
\wp(\wp(\omega)), \dots\}$, in modern notation), answering a question
by Fraenkel.  This is exactly the set the existence of which does not
follow from Zermelo's original axioms, and requires replacement, as
Fraenkel had shown.  However, the notion of ``function'' in
\citet{Fraenkel1926} does not suffice for this purpose, as
\citet{vonNeumann1928a} shows. He therefore introduced an ``extension
of the notion of function,'' the necessity of which, he says, became
clear to him in correspondence with Fraenkel (see
\citealt[p.~376]{vonNeumann1928a}). This letter is most likely part of
this correspondence, as the proof indicates how the extended notion of
function suffices for the existence
of~$\mathsf{M}$.\footnote{\Citet{vonNeumann1928a} was received by the
publisher three months prior. We may thus suppose that Fraenkel
reviewed the paper for the \emph{Mathematische Annalen}, or that von
Neumann sent Fraenkel the paper for comments. Thanks to Aki Kanamori
for clarification of the context of this letter.} See also Fraenkel's
reply to von Neumann's article \citep{Fraenkel1928a}.

The next three letters are from 1935 and 1936. They concern the
publication of \cite{Fraenkel1937}. Fraenkel had apparently asked von
Neumann if it could be published in the \emph{Annals of Mathematics},
to which von Neumann agreed in the letter of 5~December 1936. The
manuscript wasn't finished yet, as in his reply, von Neumann asked
when Fraenkel would be able to send it.  In the letter of 9~May 1936
(sent while von Neumann was travelling across the Atlantic), he
acknowledges receipt of the manuscript and proposes an extension of a
result in it. He also tells Fraenkel about the new \emph{Journal of
Symbolic Logic}, to be edited by Alonzo Church. This new journal
should become a first-rate venue for work in mathematical foundations,
and he wonders whether Fraenkel would agree to have his article
published there rather than in the \emph{Annals}. Fraenkel agreed, and
in the final reply of late May 1936, von Neumann indicates that Church
will be especially grateful and hopes to meet Fraenkel again soon,
perhaps in Jerusalem.

In the final letter of 13 November 1951, von Neumann (or presumably
his secretary, as the letter is in English) indicates that he would be
available to meet in Princeton on November 28 or~29.

\appendix
\section{Appendix: transcriptions of the letters}

In what follows, we present transcriptions of the surviving letters
from Gödel and von Neumann to Fraenkel. Editorial interpolations and
insertions are indicated by \inter{text}. Page breaks and pagination
is indicated by~\pag{n}. Text that the author has deleted (struck
through) is indicated by \sout{text}. Underlining is rendered in
\emph{italics}. Footnotes indicated by asterisks in the original are
here numbered consecutively. Handwritten letters are indicated by
``ms,'' typed letters by ``ts.'' With few exceptions, we have not
corrected von Neumann's idiosyncratic, somewhat archaic spelling.

\subsection{Fragment of a letter from Gödel to Fraenkel, c. late 1931, ms}

\selectlanguage{german}

\inter{The first two pages are missing; the first extant page is labeled as p.~3.}

\medskip

\noindent Indem man derartige Sätze mehrmals hintereinander verwendet,
kann man in der Reihe der Ord\inter{inal}-Zahlen immer höher
hinaufgelangen.--- \sout{Eines} Als Hauptbeweismittel wird der Begriff
der "`Typenmenge"'\footnote{in der üblichen Ausdruckweise = Menge der
kritischen Stellen einer Normalfunktion} verwendet. Darunter wird
verstanden die Menge derjenigen Elemente einer Menge~$M$ von
Ord\inter{inal}-Zahlen, welche gleich sind dem Ordungstypus der in~$M$
vorangehenden Elemente. Die Reihe der successiven "`Typenmengen"' kann
natürlich durch Durchschnittsbildung in's Transfinite fortgesetzt
warden. Fehler in den Beweisen sind mir nicht aufgefallen, doch wird
vielfach Bekanntes bewiesen.

Die Arbeit von Herbrand in letzten Crelle-Heft
\inter{\citealt{Herbrand1931}} habe ich \sout{bewiesen} gelesen und
stimme in allen wesentlichen Punkten (d.h. in allem, was bewiesen u.
nicht bloß vermutet wird) mit ihm überein. Der
Widerspr\inter{uchsfreiheits}-Beweis für das arithmetische System
in~§3 steht ja im wesentlichen schon in den "`Thèses"'
\inter{\citealt{Herbrand1930}}. Dieses System enthält übrigens
keineswegs alle Schlußweisen der \pag{4} klassischen Arithmetik oder
des Weylschen "`Kontinuum"' \inter{\citealt{Weyl1918}}, denn die
vollständige Induktion darf darin nur auf rekursiv definierte (daher
entscheidbare) Eigenschaften angewendet werden.\footnote{Vgl. Seite~5,
§3 (gesperrter Satz)} Daß der Herbrandsche
Wid\inter{erspruchsfreiheits}-Beweis nicht mit meinem Resultat in
Widerspruch steht hat ja Herbrand selbst (Seite 7 u.~8)
auseinandergesetzt. Er ist eben \emph{niemals} in demselben System,
dessen Widerspruchsfreiheit er beweist, formaliserbar [auch nicht in
dem Fall, wo man ihn auf das \emph{ganze} in §2 beschriebene System
anwendet], obwohl er intuitionistisch einwandfrei ist. \sout{Der}
Hilbert u. seinen Schülern kommt es natürlich gar nicht auf
Formalisierbarkeit ihrer \sout{Beweise}
Widerspr\inter{uchsfreiheits}-Beweise in irgendwelchen Systemen an,
sondern nur darauf, daß sie "`finit"' sind, wobei der Terminus
"`finit"' wohl anschaulich einen ziemlich deutlichen Sinn hat, aber
keine präzise Definition. Eine solche Definition ist sogar prinzipiell
\pag{5} unmöglich, worin ich mit den Intuitionisten und Herbrand (vgl.
S.~8) völlig einer Meinung bin. Trotz dieser zugegebenen
"`Uferlosigkeit"' der finiten Beweise \sout{glaubt} vermutet Herbrand
(und soviel ich weiß auch v.~Neumann), daß \emph{alle} finiten Beweise
schon im Axiomsystem der klassischen Analysis formalisiert werden
können, und halten daher auf Grund meines Resultats die
Undurchführbarkeit des Hilbertschen Programms für bewiesen. Mir
scheint dies nur insofern richtig, als \emph{die} Beweismittel, mit
welchen die formalistische Schule bisher den
Widerspruchfreiheitsbeweis zu führen suchte, für diesen Zweck sicher
\emph{nicht} ausreichen; im übrigen scheint mir die Herbrandsche
Vermutung (die er sogar als neues logisches Postulat vorschlägt)
durchaus unbegründet und schon ihren Sinn halte ich wegen des Fehlens
einer exakten Definition für den Be\pag{6}griff des finiten Beweises
für einigermaßen problematisch (Nebenbei: Man kann nicht etwa das
arithm. Axiomsystem in~§2 als eine Definition für den Begriff
"`finiter Beweis"' ansehen, denn der zu definierende Terminus wird ja
bei Aufstellung der Axiome selbst wieder verwendet---in Gruppe C
und~D).

Sollten Sie, sehr geehrter Herr Professor, mit Herrn v.~Neumann über
die Frage der Wid-Beweise korrespondiert haben, so wäre ich Ihnen für
eine gelegentliche Mitteilung seiner Ansicht sehr dankbar.

\closing{Mit vorzüglicher Hochachtung

\qquad Ihr ergebener\qquad Kurt Gödel}

\return{Wien VIII Josefstädterstr.~43}

\subsection{Fragment of a letter from Gödel to Fraenkel, c.~1932, ms}

\inter{The fragment consists of two sides of one sheet of paper. Only the middle of the sheet remains; the top and bottom 5 cm (approx.) are cut off.The first side contains:}

\medskip

\noindent\inter{\dots} ich Ihrem Wunsche entsprechend nach 14 Tagen
zurücksenden und erwarte mir bis dahin die Angabe Ihrer Adresse.

Zwischen meiner und der Presburgerschen Arbeit "`Über die
Vollständigkeit eines gewissen Systems\dots"'
\inter{\citealt{Presburger1930}} besteht nicht der geringste
Gegensatz. Denn in dem Presb. System ist Addition der einzige arithm.
Grundbegriff und es läßt sich daraus innerhalb des Systems die
\inter{\dots}

\medskip

\noindent\inter{The first paragraph above is struck out diagonally in
red pencil; the second one in is struck out in black pencil, and
highlighted by a vertical line in red pencil in the left margin.

The second side reads:}

\medskip

\noindent\inter{\dots Inkong?}gruenzen angeben kann---Bei
dieser Gelegenheit möchte ich noch bemerken, daß Zermelo im
Jahresbericht Bd~41, 2.~Abt. S.~87 \inter{\citealt{Zermelo1932}}
völlig unzutreffende Behauptungen über meine Arbeit aufgestellt hat,
worauf ich noch im Jahresbericht zurückzukommen hoffe.

\closing{Nochmals vielen Dank und die besten Grüße

\qquad Ihr ergebener \qquad Kurt Gödel}

\medskip

\noindent\inter{Fraenkel placed an exclamation mark in red pencil in
the right margin, next to ``völlig unzutreffende.'}

\subsection{Gödel to Fraenkel, 26 January 1958, ts}

\return{Princeton, Jan. 26, 1958.}

\smallskip

\opening{Sehr geehrter Herr Kollege:}
In Beantwortung Ihrer beiden Briefe vom 12. Dez. und 2.~Jan.
möchte ich Ihnen folgendes mitteilen:

Dr. Michael Rabin ist mir durch seine Dissertation und mehrere
Gespräche, die ich mit ihm über logische Frage hatte, bekannt. Das
genügt natürlich nicht, um mir in allen Punkten eine Meinung über
seine Eignung für eine Stellung als Assistent \inter{sic} Professor zu
bilden. Aber ich kann jedenfalls das folgende sagen:

Dr. Rabin ist zweifellos ein ausserordentlich befähigter
Logiker. Seine Dissertation enthält sehr interessante Ergebnisse.
Sie beweist ein allgemeines Theorem von grossem Interesse auf eine
überraschend einfache und höchst elegante Weise. Uber seine
Befähigung als Lehrer kann ich nichts aus eigener Anschauung sagen.
Ich glaube aber aus seiner analytischen Begabung und der grossen
Klarheit seines Denkens schliessen zu können, dass er ein
ausgezeichneter Lehrer ist. Für seine Eignung zur Zusammenarbeit
mit andern liegen ebenfalls klare Beweise vor.

Was Herrn Azriel Levy betrifft, so habe ich aus den mir zugesandten
Manuskripten ebenfalls einen sehr günstigen Eindruck von seinem
Fähigkeiten gewonnen. Bezüglich einer Stellung am Institut für 1958/59
sind die Aussichten leider sehr ungünstig, da über Stipendien an
Ausländer im wesentlichen bereits vor Weihnachten entschieden wurde
und unter den bisher nicht befriedigten Antragstellern sich sogar
solche mit wichtigen originelle Leistungen befinden. Ich werde mir
aber jedenfalls Herrn Levy's Namen vormerken für den, allerdings
unwahrscheinlichen, Fall, dass sich im Laufe der \pag{2} nächsten
Wochen doch irgend eine Möglichkeit ergeben sollte.

\closing{Mit besten Empfehlungen

\qquad Ihr

\qquad\qquad Kurt Gödel}

\subsection{Gödel to Fraenkel, 12 May 1958, ms}

\return{Princeton, 12./V.1958}

\opening{Sehr geehrter Herr Professor Fraenkel!}
Es tut mir leid, dass es heuer nicht möglich war, rechzeitig ein
Stipendium des Instituts für Herrn Levy zu erwirken. Es sind vor
kurzem unerwarteter Weise Geldbeträge verfügbar geworden, die für
Herrn Levy reserviert wurden, aber Prof. Dvoretzky teilte mir mit,
dass Herr Levy bereits andere Verfügungen getroffen hat u. daher nicht
in der Lage wäre, eine Stellung am Institut for Advanced Study
anzunehmen.

Ich hoffe, dass es möglich sein wird, das Stipendium für 1959/1960 zu
erneuen; es ist jedoch dazu eine neue Entscheidung der Fakultät
erforderlich.

\closing{Mit besten Grüssen

\qquad Ihr \qquad Kurt Gödel}

\subsection{Gödel to Fraenkel, 21 November 1958, ms}

\inter{On letterhead: Institute for Advanced Study\\
Princeton, New Jersey\\
School of Mathematics}

\return{Princeton, 21./XI. 1958}

\smallskip

\opening{Lieber Professor Fränkel!}
Besten Dank für Ihren Brief. Ich würde mich sehr freuen, Sie in
Princeton zu sehen. Da ich fast nie beim Institutstee anwesend bin,
wäre es vielleicht am besten, wenn wir die Zeit telephonisch
vereinbaren. Die Tel. No. meiner Wohnung ist: WA 4--0569

\closing{Mit besten Grüssen

\qquad Ihr \qquad Kurt Gödel}

\subsection{von Neumann to Fraenkel, 26 October 1923, ms}

\return{Zürich, den 26. X. 1923}

\smallskip

\opening{Sehr geehrter Herr Professor!}
Ich habe Ihren w\inter{ehrten} Brief und Ihre Arbeit soeben erhalten.
Entschuldigen Sie, daß Ich erst jetzt antworte, aber dieselben mußten
mir erst von Budapest hierher nachgeschickt werden.--- Ich bin Ihnen
für Ihr Interesse an der Sache und ihrer Kritik sehr dankbar, und
möchte jetzt meinen Standpunkt auseinandersetzen. Das kann ich nach
Ihrer Antwort genauer und vollständiger tun, als es in meinem ersten
Briefe möglich war.

\hbreak

Ich glaube, daß ich nicht übertreibe, wenn ich behaupte, daß
der wesentliche Teil der ganzen Mengenlehre, an dem sich die
Brauchbarkeit einer Methode zeigt, die Theorie der Ordnungszahlen ist.
(Mitinbegriffen die Alephs.) Man verlangt wohl in erster Reihe von
einer (formalistisch-axiomatischen) Mengenlehre, daß sie 
\begin{enumerate}
    \item [1.,] aus allen gleichmächtigen (bzw.
    ähnlich-wohlgeordneten) Mengen eine representative auswählen soll,
    und zwar auf axiomatisch einwandfreiem Wege;
    \item [2.,] daß die "`Vergleichbarkeit"' nachgewiesen werde.
\end{enumerate}
Daß ist aber eben die Leistung einer "`Ordnungszahlen Theorie."' Nun
ist aber die Zermelosche Axiomatik zur Herstellung des
Ordnungszahlenbegriffes unbrauchbar. Denn nichts garantirt in ihr, daß
die "`Ordnungszahl"' einer "`nicht zu großen"' Menge nicht "`zu groß"'
wird. Man weiß ja von dieser Ordnungszahl eigentlich nur, daß sie
einer gegebenen Menge ähnlich (also gleichmächtig) ist. Hieraus kann
man aber nur dann folgern, daß \pag{2} die Ordnungszahl nicht zu groß
ist, wenn man im Besitze Ihres "`Ersetzungs-Axioms"' ist. Ich war
infolgedessen stets davon überzeugt, daß dasselbe für die
formalistische Mengenlehre unbedingt notwendig ist.\footnote{Um so
mehr überrascht und interessirt es mich, daß Sie es für provisorisch
halten. Ich wäre Ihnen sehr dankbar, wenn Sie mir mitteilen wollen,
wie Sie es zu ersetzen gedenken? Daß man etwa noch die Existenz des
$\aleph_\alpha$ für jedes~$\alpha$ auf anderem Wege auch sichern
könnte, kann ich mir noch vorstellen. Es würde mich aber sehr
interessiren, ob auch der Begriff der Ordnungszahl vom Ers. Ax. frei
zu machen ist? Ich glaube kaum, daß das möglich ist. Höchstens wenn
Sie die Möglichkeit der allgemeinen Definition durch
\emph{transfinite} Induction \inter{sic} postuliren.}---

Ferner glaube ich, daß das Ers. Ax. ziemlich unbedenklich ist.
Sicherlich ist es unbedenklich in der Richtung der Richard'schen
Antinomie.--- Natürlich wäre es falsch, ohne Rücksicht auf die
Definitität, \sout{nur} es nur auf die Mächtigkeit (die tatsächliche, nicht
die axiomatische) ankommen zu lassen. Da letzten Endes doch alles
abzählbar ist, wäre das widerspruchsvoll. (D.h., man würde die
Richard'sche Antinomie bekommen.)

Aber ich fordere (wie es auch Skolem betonte) viel weniger: die
"`Ersetzung"' der Elemente einer Menge muß auf definite Art und Weise
(durch eine Operation $y = [cx]$, $c$ ein II Ding) erfolgen. Es kommt
also nicht nur roh auf den Umfang an: der definite Character bleibt
unbedingt gewahrt.---

In der Zermeloschen Axiomatik ist es jedoch kaum möglich das Ers. Ax.
unter zu bringen. Den\inter{n} die Definition der "`Function"' bei
Zermelo (Menge aller Paare $x, f(x)$) setzt ab ovo voraus, daß sowohl
der Argument- wie der Wertevorrat "`nicht zu groß"' ist. Das \pag{3} Ers.
Ax., welches fordert, daß aus dem ersteren das letztere folge, wäre
\emph{da} in der Tat überflüssig, und in seiner vollen Schärfe
garnicht ausdrückbar.

So wird man gezwungen, die Function mit \emph{beliebigem} Wertevorrat
auf andere Art, oder unabhängig einzuführen. Und ich glaube, daß es
dann am allereinfachsten ist, sofort von der "`Function"' mit
beliebigem Argument- \emph{und} Wertevorrat aus zu gehen.

Davon abgesehen kann ich zur Motivirung de\inter{s} primären
Characters der Function anführen, daß man beim formulieren des
"`Aussonderungs-Axioms"' doch auf die Function (mit beliebigem
Wertevorrat) angewiesen ist. Durch das Ausgehen aus der Function,
erspart man sich hier eine weitgehende Complication. 

\hbreak

Ich habe versucht die Gründe auseinander zu setzen, die mich veranlaßt
haben, das Ers. Ax. in die Axiomatik (explicit oder implicit) auf zu
nehmen. Das ist freilich noch keine Motivirung für das Axiom IV\,2.
Man könnte statt IV\,2 einfach das Ers. Ax. fordern, d.h. (mit meinen
Bezeichnungen):
\begin{axiom}{III\,2$^*$}
    \item[IV\,2$^*$] $a, c$ seien II Dinge, $b$ ein I~II Ding. Es
    gebe zu jedem $y$ mit $[ay] \neq A$ ein $x$ mit $[bx] \neq A$,
    $[cx] = y$.\\
    Dann ist auch $a$ ein I~II Ding.
\end{axiom}
Ich gebe zu, daß mit diesem Satz alles erreicht werden kann, was
die Mengenlehre braucht, bis auf die Auswahl. Dieselbe ließe sich,
ziemlich einfach, erreichen, wenn man III\,2, 3 durch 
\begin{axiom}{III\,2$^*$}
    \item[III\,2$^*$] $a$ sei ein II Ding. Dann gibt es ein II
    Ding~$b$, so daß man immer, wenn ein $y$ mit $[a\langle xy\rangle]
    \neq A$ existiert, $[bx]$ ein solches~$y$ ist. \pag{4}
\end{axiom}
ersetzte. Ich gebe zu, daß mit diesem System wirklich alles
erreicht werden kann, was in Frage kommt (während IV\,2 daraus nicht
folgt). Und dabei würde der Auswahl ihre Sonderstellung erhalten
bleiben.---

Trotzdem glaube ich einen Schritt über das unumgänglich notwendige
hinaus tun zu dürfen. Der Grund ist allerdings sehr subiectiv: mir
imponirt die gewaltige Leistungsfähigkeit dieses Axioms. Denn
\begin{enumerate}
    \item[1.,] Es sagt klipp und klar aus, wann eine Menge "`zu groß"'
    wird.
    \item[2.,] Es liefert ohne weiteres die Auswahl.
    \item[3.,] Und schließlich schafft es eine\inter{n} (so glaube
    ich) großen Schönheitsfehler der Mengenlehre weg: daß nähmlich zur
    Begründung der Theorie der \emph{endlichen} Mengen gewöhnlich die
    spezifisch "`\emph{unendlichen}"' Axiome~IV\,2,\,3
    \inter{axioms~V\,2, 3 are evidently meant here; see below}
    herangezogen werden müßen. (Um nämlich zu beweisen, daß Mengen wie
    $(a)$, $(a,b)$, $(a,b,c)$ nicht "`zu groß"' sind.)
\end{enumerate}

\hbreak

Man kann dieses Axiom IV\,2 auch so darstellen:

Ersetzen wir IV\,2 durch das mindestmögliche an dieser Stelle, das
Aussonderungs Axiom. Es würde so lauten:
\begin{axiom}{IV\,2$^{**}$}
    \item[IV\,2$^{**}$] $a$ sei ein II Ding, $b$ ein I~II Ding. Aus
    $[ax] \neq A$ folge $[bx] \neq A$. Dann ist auch $a$ ein I~II Ding.
\end{axiom}
Man kann dann den folgenden Satz beweisen (genau so wie die
Russel\inter{l}'{} sche Antinomie):
\begin{quote}
    Ein II Ding $a$ ist sicherlich kein I~II Ding, wenn es \pag{5} ein II
    Ding gibt, so daß zu jedem~$y$ ein $x$ mit $[ax]\neq A$, $[bx] =
    y$ existirt.
\end{quote}
IV\,2 fordert nun, daß auch die Umkehrung richtig sei, daß diese
hinreichende Bedingung auch notwendig sei. D.h.: daß \emph{möglichst
viele} Bereiche Mengen seien.---

IV\,2$^*$ (das Ers. Ax.) kann man folgendermaßen zu IV\,2 ergänzen
(wobei die Auswahl, d.h. III\,2$^*$ nicht vorausgesetzt wird):
\begin{quote}
    Es soll eine Wohlordnung $W$ des Bereichs der I Dinge~$\Omega$
    geben, bei der \emph{alle} echten Abschnitte von~$\Omega$ Mengen
    (nicht nur Bereiche) sind.
\end{quote}
Dies folgt aus IV\,2 (auf Seite 39, oben: man wähle $W$ so, daß
$\Omega, W \approx \Omega^{**}, \overline{\Sigma}(\Omega^{**})$ wird).
Aber umgekehrt ergibt das mit IV\,2$^*$ das Axiom~IV\,2.

\hbreak

Ich glaube daß die Axiomen Gruppe V (Unendlichkeits Axiome) ein großer
Schönheitsfehler des ganzen Systems ist. Sie sind aus dem Zermeloschen
System übernommen und sehr complicirt. Ihre Bedeutung (besonders
V\,2\,3) ist mir absolut unkar. Sie sind aber notwendig, um
\begin{enumerate}
    \item [1.,] die Existenz der unendlichen Ordnungszahl~$\omega$ zu zeigen,
    \item [2.,] um zu zeigen, daß der Bereich aller Alephs,~$\Gamma$,
    keine Menge, also "`zu groß"' ist.
\end{enumerate}
Dies auf andere Art zu sichern, ist mir noch nicht gelungen. Es wäre
aber sehr erwünscht, denn vielleicht ist die Quelle der
Schwierigkeiten des \emph{Continuum-Problems} hier.

\hbreak

Der Grund, die Aequivalenz \emph{nach} der Ähnlichkeit zu
behan\pag{6}deln, ist der: man braucht schon für die einfachsten Sätze
der Aequivalenz (Definition der Alephs, Vergleichbarkeit,
Aleph-Operationen) die ganze Theorie der Ordnungszahlen und den
Wohlordnungssatz.

Umgekehrt setzt die Theorie der Wohlordnung nichts über die
Aequivalenz voraus. Es wäre unschön das ganze Kapitel~II zwischen den
Bernsteinschen Satz und die Vergleichbarkeit einzuschieben. (Bei der
gegenwärtigen Anordnung geht man immer mehr und mehr vom Allgemeinen
aufs Spezielle über.)

Der Begriff der Ordnungszahl ist eben (wenigstens in meiner
Darstellung) das wesentliche, alles andere ist secundär.

\hbreak

Ich möchte die Gelegenheit benützen, um einen Irrtum in meinem Aufsatz
zu corrigieren. Auf Seite~39 unten behaupte ich, daß die
Vergleichbarkeit von wohlgeordneten Bereichen sich auch dann
nachweisen läßt, wenn dieselben keine Mengen mehr sich ("`zu groß"'
sind). Der Beweis wäre ohne Ordnungszahlen zu führen, etwa indem man
je zwei Elemente mit ähnlichen Abschnitten einander zuordnet.

Das ist falsch. Alle Methoden versagen daran, daß auch die echten
Abschnitte keine Mengen mehr zu sein brauchen. (Es ist unmöglich die
Existenz eines II~Ding~$c$ das die Abschnitte von $x$ und~$y$
aufeinander ähnlich abbildet als Function von $x, y$ aus zu drücken.
Denn $c$ braucht kein I~Ding zu sein.)

Das ist wohl ein weiterer Beleg dafür, daß nicht die Wohlordnung,
sondern die Ordnungszahl der zentrale Begriff ist: \pag{7}
wenn die Ordnungszahl fehlt, hat der wohlgeordnete Bereich gar nichts
mehr vor den übrigen voraus.

\hbreak

In dem Teile II  fehlt der Beweis eines Satzes, der auch hierher
gehört, um so mehr als er vielfach für selbstverständlich gehalten
wird. Es ist die Möglichkeit der Definition durch transfinite
Induction. In der weitesten Fassung lautet er so:
\begin{quote}
    $a$ sei ein II Ding. Es gibt dann ein II Ding~$b$, und nur eins,
    mit den folgenden Eigenschaften:
    \begin{enumerate}
        \item [1.,] Es ist $[bx] = A$, wenn nicht $x$
        O\inter{rdnungs}Z\inter{ahl} ist.
        \item [2.,] Es ist $[bx] = [a\langle
        x\left|[bx]\right|\rangle]$, wenn $x$
        O\inter{rdnungs}Z\inter{ahl} ist.
    \end{enumerate}
\end{quote}
Der Beweis ist dem Beweis der Existenz der Ordnungszahlen weitgehend
analog. (Dort ist $a$ sehr einfach gewählt: es ist stets $[a\langle
xy\rangle] = y$. Die Allgemeinheit $a$-s muß hier berücksichtigt
werden. Vereinfachungen hat man hierbei, weil der Begriff der
Ordnungszahl schon vorhanden ist.)

\hbreak

Das ist es ungefähr, was ich über die Teile I, II sagen kann. Ich will
nun auf den Teil~III übergehen.---

Das Verhähltnis dieser axiomatischen Mengentheorie zum
"`Intuitionismus"' wird in diesem Teil~III so gedacht:

Wir befinden uns in einem Bereiche (Bereich soll von nun an nicht mehr
den \sout{Sinn} axiomatischen Sinn habe, es ist einfach der naïve
Collectivbegriff), dem die I und II~Dinge angehören. In ihm sind die
elementaren Operationen der Logik:
\begin{quote}
    logisches Schließen, Negation, Gleichheit, Begriffe "`alle"'
    \pag{8} und "`es gibt"'
\end{quote}
schon sinnvoll festgelegt, und diese werden auch nicht mehr analysirt
werden; demgegenüber fehlen noch vollständig die der
"`Imprädikativität"' verdächtigen Begriffe der
\begin{quote}
    Menge, Function, Auswahl.
\end{quote}
Diese sollen axiomatisirt werden.

\hbreak

Dieser Bereich ist also nicht ganz "`intuitionistisch"', wie etwa
diejenigen aus denen Brouwer oder Hilbert vorgehen. Er entspricht
vielmehr dem ersten, weniger radikalen, Standpunkte Weyl's, und ist
von ihm in seinem Buche "`Das Kontinuum"' wohl zuerst \emph{genau}
beschrieben worden.---

\hbreak

Ich weiß, daß das ein halber und überholter Standpunkt ist: um die
Mengenlehre zu reconstruiren, muß man sich, wie Hilbert, rückhaltlost
auf den Boden des intransingentesten Brouwerschen Intuitionismus
stellen. (Solange man "`inhaltlich schließt"'.) Da es mir aber nicht
auf den Nachweis der Widerspruchsfreiheit ankommt, sondern nur darauf,
überhaupt anzugeben, was die formalistische Mengenlehre ist, tue ich
es lieber nicht.

Die Axiome lassen sich in einem solchen halb-\-intu\-itionist\-ischen
Bereich (wie i\inter{h}n auch Zermelo benützt hat) viel leichter
formuliren. Außerdem kann man \sout{dieses System} diesen Bereich
(Logik ohne Collectivbegriffe) leicht im Sinne des Hilbertschen
Verfahrens herstellen, und seine Widerspruchsfreiheit
nach\pag{9}weisen. (Hilbert, Logische Grundlagen der Mathematik, Math.
\sout{Zeits.} Annalen Bd.~88 \inter{\citealt{Hilbert1923}}. Hier
kommen in Frage: Axiome I\,1--4, II\,5--6, III\,7--8, V\,11.) Man hat
dann die wesentliche Frage, das Mengen\-lehren-Ax\-iomen\-system
isolirt.---

Man kann übrigens auch das ganze direct in das Hilbertsche System
eingliedern. Einige \sout{kleine} Abweichungen von der ursprünglichen
Form treten dabei auf.

\hbreak

Die Abzählbarkeit des Systems steht zur Existenz von (axiomatisch)
unabzählbaren Mengen in keinem Widerspruch. Wie es auch Skolem bemerkt
hat, bedeutet die axiomatische Unabzählbarkeit gar nicht die
tatsächliche Unabzählbarkeit.---

Ich möchte hier das folgende bemerken: beim Abfassen des Aufsatzes war
mir die Skolemsche Arbeit noch nicht bekannt. Ich habe eine
überflüssige Complication beim Beweise der \sout{Un} Abzählbarkeit
eingeführt (beim Axiom V\,1): darum konnte ich bloß nachweisen, daß
der Bereich (naïv) nur der I und~II Dinge dem $\omega$ (tatsächlich)
aequivalent ist. Nach dem ich die Skolemsche Arbeit eingesehen hatte
(ich bin Herrn Professor sehr dankbar dafür, mich auf dieselbe
verwiesen zu haben), constatirte ich dies. Man kann, wie ich nachher
sah, mit denselben Mitteln auch \emph{direct} die Aequivalenz mit der
"`\emph{Folge}"' $0, F(0), F(F(0)), F(F(F(0))), \dots$ (also $1, 2, 3, \dots$;
das ist ein\sout{e} \emph{Teil}\sout{menge} von~$\omega$, eventuell
ein\sout{e} echter) nachweisen. Man behandelt nur dazu das Axiom~V\,1
ebenso wie die übrigen Axiome (etwa wie III\,2, 3 oder IV\,2).

\hbreak

Bezüglich der "`Relativität"', oder wie ich sagte
"`Scheinbar\pag{10}keit"' der Mächtigkeiten stehe ich ganz auf Skolems
Standpunkt. Allerdings ohne dieselbe Consequenzen hieraus zu ziehen
wie er: ich glaube vielmehr, daß das ein schweres Argument
\emph{gegen} den Intuitionismus ist. Nie wird man auf
Intuitionistischer Basis dem Begriff des unabzählbaren Continuums auch
nur nahe kommen können, denn die wirkliche Ideenwelt ist abzählbar.
Das Continuum ist in seinem Wesen transfinit, imprädikativ, und
\emph{nur} für die formalistische Methode faßbar.

\hbreak

Ich bin, wie aus dem obigen hervorgeht bezüglich der Imprädikativität
alles unabzählbaren vollkommen Ihrer Meinung. Man kommt mit der
intuitionistischen Methode nicht durch: man muß sozusagen den
formalistischen Apparat (wie ihn Hilbert andeutet) zwischen das
unmittelbare ("`inhaltliche"', intuitionistische) Schließen und die
Mengenlehre bzw. Mathematik einschieben.---

Und das ist nur zulässig, wenn der Apparat widerspruchsfrei ist. Es
ist nun leicht nach zu weisen, daß dies der Fall ist, solange \inter{n}ur 
\begin{quote}
    das Folgern,\footnote{Das Folgern, $a \to b$, ist vom Schließen,
    $\begin{array}{c}a\\ a \to b\\ \hline b\end{array}$, streng zu
    unterscheiden!} die Negation, die Identität, die
    Auswahl\footnote{Die ganz allgemeine Auswahl ist auch
    unbedenklich, solange das wa\inter{h}re Imprädikative Element, die
    \emph{Function}, fehlt. Man kann die Widerspruchsfreiheit solange
    leicht nachweisen.}
\end{quote}
darin vorkommt.

Will man aber die Mengenlehre erhalten, so muß man die
Function\footnote{Hilbert verwirrt den Tatbestand ziemlich, indem er
einen neuen Begriff, die Prädikate~$A(a)$, einführt. Man vermeidet
dieselben besser, in dem man bei den gewöhnlichen Ausdrücken~$a$
bleibt, und die Substitution einer Variabeln, $\mathrm{Subst}({b
\atop x})a$, definirt. Etwas Vorsicht ist geboten, eine Variable kann,
z.B. wegen der Auswahl ev. nicht "`frei"' sein.)\par Die "`Function"'
wären zwei Operationen $f, \phi$ derart, daß $\mathrm{Subst}({b
\atop x})a = \phi(f_{a,x},b)$ Axiom ist. (Allgemein führt das
freilich zur Russel\inter{l}schen Antinomie. Man muß eben wieder
einschränken.)} \pag{11} mit den entsprechenden Axiomen einführen. (Meine
Axiomatik z.B. läßt sich leicht so umformen. Die Teile II, III, IV\,1
fallen automatisch weg.) Die Widerspruchslosigkeit für \emph{diesen}
Formalismus ist aber noch nicht nachgewiesen.---

Ich wei\inter{ß} nicht einmal, ob Hilbert diesen Weg eingeschlagen hat.
Vielleicht will er die Mengenlehre "`stufenweise"' herstellen. Ich
glaube daß diese "`Umgehung der Imprädikabilität"' kaum möglich sein
wird. Russel\inter{l}s Fiasko beweist es.

Allenfalls fehlt hier, wo die formalistische Mengenlehre anhebt, ein
Widerspruchsfreiheitsbeweis. Diese Lücke ist aber wohl das, was man
als "`Imprädikabilität"' zu bezeichnen pflegt.

\hbreak

Zum Schluße möchte ich noch drei Bemerkungen machen.---

Das Axiomensystem ist \emph{nicht} kategorisch. Das ist ja klar, ich
fordere nicht einmal, daß alles Function sei. Noch weniger sind Mengen
mit Eigenschaften wie $a = (a)$, $a = (((a)))$, $a=(a ((a)))$, etc.
ausgeschlossen. (Das sind wohl die "`Ensembles extraordinaires"' des
Herrn Mirimanoff?) Allein \emph{nicht das} ist das Wesentliche.

Denn man kann leicht das System so einengen, daß \sout{jede}
\begin{enumerate}
    \item [1.,] alle I Dinge auch I~II Dinge sind,
    \item [2.,] keine "`Absteigende Functionen-Folge"' existirt (wie
    bei Skolem), d.h.: Es gibt kein II~Ding~$a$, so daß für jede
    N\inter{aturliche}Z\inter{ahl}~$p$ \quad $[a F(p)]$ einem~$x$ oder
    einem $[[ap]x]$ mit $[[ap]x] \neq A$ gleich ist.\footnote{Bei
    Skolem würde es lauten: \dots, so daß für jede NZ~$p$ $[a F(p)]$
    ein Element von $[ap]$ ist. Also: daß $[a F(p)]$ gleich einem~$x$
    mit $[[ap]x] \neq A$ ist.---Das wäre bei \emph{Mengen} so. Bei
    Functionen muß ausser~$x$ auch $[bx]$ in Betracht genommen
    werden.---\par Der Satz gilt nur mit gewissen Beschränkungen, vgl.
    folgende Note.}~\pag{12}
\end{enumerate}
Hierdurch hätte man erreicht, daß jede Function in einer
wohlgeordneten Reihe von Schritten aus $0$ und $(0)$\footnote{Man
könnte noch leicht $A =0$, $B=(0)$ erreichen. Hierauf will ich nicht
eingehen. Dann währen $0$, $(0)$ characterisirt durch das folgende:
\begin{quote}
    $0$ is ein $a$, so daß stets $[ax]=a$ ist.\\
    $(0)$ is ein $b$, so daß $[ba]=b$, und sonst $[bx] = a$ ist.
\end{quote}
Diese beiden sind also auf nichts früheres zurückführbar.} aufgebaut
werden kann. Die vorhin angeführten barokken Mengenbildungen würden
wegfallen.---

Trotzdem hilft das alles nichts; das System ist noch immer nicht
kategorisch. Der Grund mag in dem "`relativen"' Character der
Wohlordnung liegen.

Nach alledem glaube ich, daß sich die Mengenlehre gar nicht
kategorisch axiomatisiren läßt. Dann gibt es aber (von endlichen
Systemen abgesehen) überhaupt keine kategorische Axiomatik auf der
Welt. (Scheinbar kategorische, wie z.B. die Mathematik, sind von der
Mengenlehre abhängig.) Über die Tragweite und die Bedeutung dieser
Tatsache kann ich mir zunächst gar keinen Begriff machen.---

Übrigens braucht es auch keine "`kleinsten"' Systeme zu geben. D.h.:
ich kann nicht die Existenz eines Systems nachweisen, welches meinen
Axiomen genügte, während keines \sout{ihrer T} seiner Teilsysteme
diesen genügt. (Bei unverändertem Beibehalten von $A$, $B$, $[x,y]$,
$\langle x,y\rangle$, siehe Seite~48.)

\hbreak

Der Begriff der Wohlordnung, oder spezieller der der Endlichkeit ist
relativ. D.h.:

Es sei ein Complex~$\mathfrak A$ von Dingen gegeben. $\Sigma$ sei ein
den Axiomen genügendes System. Wenn nun
\begin{enumerate}
    \item [1.,] alle Elemente von $\mathfrak A$ I~Dinge in~$\Sigma$ sind,
    \item [2.,] \sout{ein II Ding~$\mathfrak A_\Sigma$} ein Bereich in $\Sigma$, $\mathfrak A_\Sigma$,
    existirt, dessen Elemente mit denen von \pag{13} $\mathfrak A$
    identisch sind,
\end{enumerate}
so heiße $\mathfrak A$ zu~$\Sigma$ gehörig.

Wenn nun $\mathfrak A$ ein beliebiges Complex ist, so kann es mehrere
Systeme $\Sigma'$, $\Sigma''$, $\Sigma'''$, \dots geben, zu
denen~$\mathfrak A$ gehört. Wenn $\mathfrak A$ endlich ist (naïv) so
sind die $\mathfrak A_{\Sigma'}$, $\mathfrak A_{\Sigma''}$, $\mathfrak
A_{\Sigma'''}$, \dots (axiomatisch in bzw. $\Sigma'$, $\Sigma''$,
$\Sigma'''$, \dots) auch endlich.---

Die Umkehrung gilt aber nicht.

So kann für dasselbe $\mathfrak A$ in einem~$\Sigma$ das $\mathfrak
A_\Sigma$ endlich sein, und in einem anderen unendlich. (Dabei muß es
tatsächlich unendlich sein.)

So ist es denkbar, daß $\mathfrak A_\Sigma$ in~$\Sigma$ zwar unendlich
ist, aber im Teil-Systeme~$\Sigma'$ von~$\Sigma$ endlich erscheint.
D.h.: beim \emph{verfeinern} des Systems~$\Sigma'$ zu~$\Sigma$ stellt
es sich von \emph{scheinbar} endlichen~$\mathfrak A$ heraus, daß
\sout{es unendlich ist} sie
unendlich sind.---

Wenn man nun nicht gerade zu intuitionistisch ist, so kann man fragen:

Gibt es ein "`ideales"' System, so daß es sich bei keiner weiteren
Verfeinerung mehr von endlichen Mengen herausstellen kann, daß sie
unendlich sind?

Oder ist jedes System, wie man es auch wählt, verfeinerungsfähig,
derart, daß man endliche Mengen als unendlich erkennt?

D.h.:

Entspricht wenigstens dem anschaulichen Begriff der Endlichkeit etwas
formales? Oder is das "`unendliche"' nach unten hin ebenso relativ und
verschwommen begrenzt, wie wir es nach oben bereits wissen?---

Ich glaube, obzwar mir der Beweis fehlt, daß das letztere \pag{14} der
Fall ist. Wenn sich das nachweisen ließe, so wäre es ein vernichtendes
Argument \emph{gegen} den vom Intuitionismus verfochtenen
anschaulichen Character der finiten Induction und der positiven
ganzen Zahl.

\hbreak

Ich habe in den bisherigen Zeilen meine Ansicht über die
formalistische Mengenlehre und ihre Axiomatisirung zu entwickeln
versucht. Ich wäre Herrn Professor sehr dankbar, wenn Sie mir Ihre
Ansicht darüber mitteilen würden.---

Die Publikation der Arbeit is mir von Herrn Professor Erhard Schmidt
im Sommer in Berlin versprochen worden. Und zwar im vollen Umfange
(ca. 40 Druckseiten) in der Math\inter{ematischen} Zeitschrift. Ehe ich
veröffentliche, möchte ich die Arbeit noch einmal umstylisiren und
einige Änderungen, die notwendig sind, vornehmen. Ich hoffe das bis
Mitte November zu erledigen. Halten Herr Professor eine vorläufige
Publication auch so für zweckmässig? Ich wäre dazu gerne bereit.
Welche Zeitschrift käme da in Betracht?---

Ich möchte Herrn Professor fragen, ob die Ordnungszahlen-Theorie, die
ich aufstelle, neu ist. Oder kannte sie schon Zermelo oder jemand
anderer früher?

\hbreak

\closing{Indem ich ihnen für Ihr Interesse nochmals danke,

\qquad
verbleibe ich hochachtungsvoll

\qquad\qquad
Ihr ergebener \qquad Johann v. Neumann}

\smallskip

\noindent Zur Zeit: Zürich, Plattenstr.~52 bei Degen. Ich bleibe hier
bis März 1924, mit Ausnahme von 20 December--5 Januar. Dann bin ich in
Budapest, Vilmos császár út 62.III

\subsection{von Neumann to Fraenkel, 17 June 1927, ms}

\inter{On Letterhead: H\^otel Bristol\\
Telegramm-Adresse: Bristol, Karlsbad}

\smallskip

\return{Karlsbad, den 17. 6. 1927}

\opening{Sehr geehrter Herr Professor!}
Ich habe Ihre Karte soeben erhalten, und beeile mich Ihre Frage zu
beantworten.

$\mathfrak{Z}$ sei die Menge der positiven ganzen Zahlen,
$\mathfrak{A}(x)$ die Potenzmenge von~$x$, $[x,y]$ das
\emph{geordnete} Paar von $x$ und~$y$. (Natürlich sind
$\mathfrak{A}(x)$ und $[x,y]$ durch Funktionen darstellbar.)

Ich nenne eine Menge~$\mathsf{M}$ einen "`Anfang"', wenn sie die folgende
Eigenschaften hat:

\begin{quote}
Jedes Element $u$ von $\mathsf{M}$ ist von der Form
\begin{align*}
    u & = [1, \mathfrak{Z}]\\
    \intertext{oder}
    u & = [x+1, \mathfrak{A}(y)]
\end{align*}
wobei aber auch $[x,y]$ zu~$\mathsf{M}$ gehört.\inter{von Neumann drew a
vertical line in the left margin along the above.}
\end{quote}

Man bringt diese Bedingung leicht auf die Form \[f(\mathsf{M}) = 0.\tag*{\mbox{\pag{2}}}\]

Man zeigt: wenn $\mathfrak{M}$ eine Mengen von "`Anfängen"' ist, so
ist der Durchschnitt von~$\mathfrak{M}$ auch ein "`Anfang"'; wenn also
$u$ zu einem "`Anfange"' gehört, so gibt es einen
"`Anfang"'~$\mathsf{M}$ mit dieser Eigenschaft:
\begin{quote}
    $u \mathrel{\varepsilon} \mathsf{M}$; für keinen
    "`Anfang"'~$\mathsf{M}'$ der echte Teilmengen von~$\mathsf{M}$
    ist, ist $u \mathrel{\varepsilon} \mathsf{M}'$.\inter{Line along
    left margin.}
\end{quote}

Auch diese Eigenschaft bringt man ohne weiteres auf die Form
\[g(u,\mathsf{M}) = 0.\] Auf Grund seiner Definition existirt für
jedes~$u$ kein oder nur ein~$\mathsf{M}$ mit dieser Eigenschaft.

Also bilden alle $\mathsf{M}$ mit \[g(u,\mathsf{M}) = 0\] eine Menge,
und \emph{hier erst} wende ich die besprochene Erweiterung des
Funktionsbegriffs an: die Menge aller dieser~$\mathsf{M}$ ist $= h(u)$.

Die Vereinigungsmenge von $h(u)$ ist~$k(u)$, es ist leer, oder der
kleinste $u$ enthaltende Anfang.~\pag{3}

\sout{Nun bildet man das durch $k$ vermittelte Bild von
$\mathfrak{Z}$} Man beweist leicht, dass für jedes $x$
von~$\mathfrak{Z}$ ein $x$~enthaltender "`Anfang"' existirt
(Induktions-\emph{Beweis}), also enthält $k(x)$ ein Element~$[x,y]$:
$y$~ist offenbar die $x$-mal iterirte Potenzmenge von~$\mathfrak{Z}$.
Man könnte es als~$j(x)$ darstellen.

Am schnellsten aber kommt man so durch: unter Anwendung des
Ersetzungs-Axioms sei~$\mathsf{Z}$ das durch $k(x)$ vermittelte Bild
von~$\mathfrak{Z}$; die Vereinigungsmenge von~$Z$
sei~$\overline{\mathsf{Z}}$. Dieses $\overline{\mathsf{Z}}$ ist, wie
man sofort zeigt, die gewünschte Menge \[\{\mathfrak{Z},
\mathfrak{A}\mathfrak{Z}, \mathfrak{A}\mathfrak{A}\mathfrak{Z},
\mathfrak{A}\mathfrak{A}\mathfrak{A}\mathfrak{Z}, \dots\}.\]

\hbreak

Im übrigen verlasse ich Göttingen am 1.7. endgültig, ich habe
einen Vortrag in Kiel am 2.7., und gehe dann nach Berlin. Um Pfingsten
habe ich in Königsberg 3 Vorträge über die
Wi\pag{4}derspruchsfreiheit gehalten.

\closing{In der Hoffnung, Sie möglichst bald wiederzusehen

\qquad verbleibe ich Ihr ganz ergebener

\qquad\qquad J. v. Neumann}

\subsection{von Neumann to Fraenkel, 5 December 1935, ms}

\inter{On letterhead: Institute for Advanced Study\\
Princeton, New Jersey}

\return{Dec. 5.}

\opening{Lieber Herr Professor Fraenkel,}
vielen Dank für Ihren lieben Brief. Wir sind selbstverständlich
hocherfreut, dass Sie daran denken, eine Arbeit in den ``Annals'' zu
publizieren, und ich kann Ihnen die gewünschte prinzipielle Zusage
im Namen der Redaktion gerne geben. Wann können sie uns das Mskr.\
schicken?

Farkas bin ich im Sommer in Budapest begegnet, und \pag{2} ich habe
von ihm wieder gehört---was ich ohnehin wusste---wie interessant und
anziehend Palästina ist. Auch ich hoffe es früher oder später besuchen
zu können.

Uns geht es gut, seit 9 Monaten bin ich Vater einer Tochter, die
Marina heisst.

Besteht keine Aussicht, Sie ausserhalb von Palästina wiederzusehen?
Gehen Sie nach Oslo?

Ich habe im Mai 1936 Vorlesungen in Paris zu halten, und bin dann
für 2--3 Sommermonate in Ungarn.

\closing{Mit den besten Grüssen,

\qquad Ihr ergebener J. v. Neumann}

\subsection{von Neumann to Fraenkel, 9 May 1936}

\inter{On letterhead: United States Lines\\
On Board S.S. ---}

\return{May 9., 1936}

\opening{Lieber Herr Professor Fraenkel,}
vielen Dank für Ihren Brief und für das Manuskript, dessen Inhalt
mich übrigens auch persönlich sehr interessiert hat.

Im Zusammenhange damit möchte ich Sie noch fragen, ob die Erledigung
des Falles mit 2-Element-Mengen nicht doch miteingeschlossen werden
könnte. Nicht-axiomatisch ist es ja ziemlich einfach:

Nehmen wir an, das Produktaxiom wäre für alle jene
Mengen~$\mathfrak{M}$ gesichert, deren Elemente paarweise
elementfremde endliche Mengen mit $\ge 3$ Elementen sind. Dann kann
man Mengen~$\mathfrak{N}$, deren Elemente paarweise elementfremde
endliche Mengen mit $\ge 1$ Elementen sind, so erledigen: \pag{2}

Seien $\overline{1}$, $\overline{2}$, $\overline{3}$ drei verschiedene
Objekte $\mathrel{\not\!\varepsilon}\mathfrak{S}(\mathfrak{N})$, sei
$f(x)$ eine Funktion derart, dass \[f((x,u))=x\] falls $x \neq
\overline{1}$, $\overline{2}$, $\overline{3}$ und $u = \overline{1}$
oder $\overline{2}$ oder~$\overline{3}$. Sei ferner $g(a)$ eine
Funktion derart, dass \sout{$g(a)$} falls $a$ eine Menge ist, $g(a)$
die Menge aller $(x,u)$, $x \mathrel\varepsilon a$, $u = \overline{1}$ oder
$\overline{2}$ oder~$\overline{3}$ ist. ($g(a) \subset a +
(\overline{1}, \overline{2}, \overline{3})$.)

Die Menge $\mathfrak{M}$ der $g(a)$, $a \mathrel\varepsilon
\mathfrak{M}$, ist leicht gebildet, \emph{ohne} Ersetzungs-Axiom, da
sie $\subset \mathfrak{P}(\mathfrak{P}(\mathfrak{S}(\mathfrak{N}) +
(\overline{1}, \overline{2}, \overline{3})))$ ist. Anwendung des
Produkt\-axioms auf~$\mathfrak{M}$ gibt ein~$\mathfrak{X}$, das mit
jedem Element von $\mathfrak{M}$ genau ein Element gemein hat. Das
$f(x)$-Bild $\mathfrak{Y}$ von~$\mathfrak{X}$ (\emph{ohne}
Ersetzungsaxiom, denn es ist $\subset \mathfrak{S}(\mathfrak{N})$) hat
daher mit jedem Element von $\mathfrak{N}$ genau ein Element gemein.

Wird diese Betrachtung in der axiomatischen Formalisierung wesentlich
schlimmer?--- \pag{3}

\smallskip
\noindent\inter{In top margin, possibly in Fraenkel's hand:}\\
1) Sonderabz.\\
2) Eingangsdatum
\smallskip

In den Monaten, die zwischen Ihrem ersten und zweiten Briefe vergangen
sind, ist bei uns die folgende Änderung eingetreten: Eine neue
Zeitschrift, ``Journal of Symbolic Logic'', wurde gegründet, deren
Redakteur A.~Church, Princeton, ist. Wir haben im Zusammenhange damit
beschlossen, alle logisch-mengentheoretische-Grundlagen Abhandlungen in
unserem Bereich nicht mehr in den ``Annals'' sondern im ``J.S.L.'' zu
drucken. In diesem Zusammenhange möchte ich Sie daher fragen, ob es
Ihnen angenehm wäre, falls auch Ihre Arbeit dort erschiene. Church
wäre sehr froh, wenn er Ihre Arbeit bekommen würde, und es ist mit
Bestimmtheit zu erwarten dass das ``J.S.L.'' eine erstklassige
Zeitschrift (übrigens wohl die erste dieser Art) sein wird.

Die erste Nummer des ``J.S.L.'' erscheint \pag{4} in einigen Wochen.
Church versprach mir, dass Ihre Arbeit in der 2-ten oder 3-ten Nummer
(3~oder 6~Monate später) erscheinen würde.

Betr.~der Korrekturen wäre es vielleicht am besten, wenn Sie mit
Church direkt in Verbindung träten. (Princeton University, Dept. of
Math., Fine Hall---Princeton, N.J., USA)---

Wie Sie sehen bin ich auf dem Wege nach Europa.  Ich habe in Paris 2
Wochen lang Vorträge, dann gehe ich vielleicht nach Ungarn. Meine
Pläne sind noch recht unbestimmt, und vielleicht kehre ich noch vor
dem Ende des Sommers nach U.S.A. zurück. Nach Oslo gehe ich
voraussichtlich nicht.

Aber ich hoffe sehr, dass wir uns in nicht allzulanger Zeit doch noch
persönlich wiedersehen.

\closing{Mit den besten Grüssen

\qquad bin ich ihr sehr ergebener \qquad J. v. Neumann}

\smallskip\noindent
\inter{In left margin:}
Meine Adresse ist bis zum 26 Mai: Institut Henry \inter{sic} Poincaré,
Paris. Nachher: Die Princetoner Adresse, von wo mir die Post
nachgeschickt wird.

\subsection{von Neumann to Fraenkel, May 1936, ms}

\inter{On letterhead: Université de Paris, Faculté des Sciences,
Institute Henri Poincaré, 11, Rue Pierre-Curie (V\textsuperscript{e}),
Tél: Odéon 42-10}

\return{Paris, le Mai 1936}

\opening{Lieber Herr Professor Fraenkel,}
herzlichsten Dank für Ihren lieben Brief, der mich noch hier
erreicht hat. Ich danke Ihnen für die liebenswürdige
Bereitwilligkeit, Ihre Abhandlung auch dem J. f. S. L. zu geben, und
Church wird Ihnen ganz besonders verbunden sein. Ich nehme an, dass
mit den Separaten alles in Ordnung sein wird. Und ich hoffe sehr,
dass Sie auch uns nicht vergessen werden, und sich noch eine
Gelegenheit ergeben wird, dass die ``Annals'' etwas von Ihnen
bekommen.---

Ich habe hier (am Inst. H. P.) Vorträge gehalten, und fahre am 4. Juni
nach Princeton zurück, während meine Frau noch einige Zeit in Europa
bleibt. Nach \pag{2} Oslo gehe ich nicht, es tut mir wirklich sehr
leid, dass damit eine Gelgenheit des Wiedersehens vorbeigeht.
Hoffentlich sehen wir uns bald anderswo wieder, vielleicht in
Jerusalem!

Mit den besten Grüssen, denen sich unbekannterweise auch meine Frau
anschliesst,

\closing{Ihr stehts ergebener 

\qquad J. v. Neumann}

\subsection{von Neumann to Fraenkel, 13 November 1951, ts}

\inter{On letterhead: Institute for Advanced Study\\
Princeton, New Jersey\\
School of Mathematics}

\return{November 13, 1951}

\opening{Dear Colleague,}
Thank you for your note of November~1. I expect to be in Princeton on
November 28 and~29, and I very much hope that we can see each other at
that time. Please let me know when you arrive here.

With best personal regards, in which Mrs.~von Neumann joins me, I am,

\closing{\qquad Cordially yours,

\qquad\qquad John von Neumann}

\smallskip
\noindent JvN:eg

\selectlanguage{english}


\begin{thebibliography}{62}
\providecommand{\natexlab}[1]{#1}
\providecommand{\url}[1]{\texttt{#1}}
\providecommand{\urlprefix}{URL }
\expandafter\ifx\csname urlstyle\endcsname\relax
  \providecommand{\doi}[1]{doi:\discretionary{}{}{}#1}\else
  \providecommand{\doi}{doi:\discretionary{}{}{}\begingroup
  \urlstyle{rm}\Url}\fi

\bibitem[{Ad{\v z}i{\'c} and Do{\v s}en(2016)}]{AdzicDosen2016}
Ad{\v z}i{\'c}, Milo{\v s} and Kosta Do{\v s}en, 2016.
\newblock G\"odel's {{Notre Dame}} course.
\newblock \emph{The Bulletin of Symbolic Logic}~22(4): 469--481.
\newblock \doi{10.1017/bsl.2016.36}.

\bibitem[{Bellotti(2016)}]{Bellotti2016}
Bellotti, Luca, 2016.
\newblock Von {{Neumann}}'s consistency proof.
\newblock \emph{The Review of Symbolic Logic}~9(3): 429--455.
\newblock \doi{10.1017/S1755020316000198}.

\bibitem[{Bhattacharya(2021)}]{Bhattacharya2021}
Bhattacharya, Ananyo, 2021.
\newblock \emph{The Man from the Future: The Visionary Life of {{John}} von
  {{Neumann}}}.
\newblock {London}: {Penguin}.

\bibitem[{Budiansky(2021)}]{Budiansky2021}
Budiansky, Stephen, 2021.
\newblock \emph{Journey to the Edge of Reason: The Life of {{Kurt G\"odel}}}.
\newblock {New York, NY}: {Norton}.

\bibitem[{Crocco and Engelen(2016)}]{CroccoEngelen2016}
Crocco, Gabriella and Eva-Maria Engelen, eds., 2016.
\newblock \emph{Kurt {{G\"odel}} Philosopher-Scientist}.
\newblock {Presses universitaires de Provence}.

\bibitem[{Dawson(1985)}]{Dawson1985}
Dawson, John~W., Jr., 1985.
\newblock Completing the {{G\"odel-Zermelo}} correspondence.
\newblock \emph{Historia Mathematica}~12(1): 66--70.
\newblock \doi{10.1016/0315-0860(85)90070-9}.

\bibitem[{Dawson(2003)}]{Dawson2003}
Dawson, John~W., Jr., 2003.
\newblock Introductory note to the {{G\"odel-Zermelo}} correspondence.
\newblock In  \parencite{Godel2003a}, 419--421.

\bibitem[{Ebbinghaus and Peckhaus(2015)}]{EbbinghausPeckhaus2015}
Ebbinghaus, Heinz~Dieter and Volker Peckhaus, 2015.
\newblock \emph{Ernst {{Zermelo}}}.
\newblock {Berlin}: {Springer}.
\newblock \doi{10.1007/978-3-662-47997-1}.

\bibitem[{Ewald(1996)}]{Ewald1996}
Ewald, William~Bragg, ed., 1996.
\newblock \emph{From {{Kant}} to {{Hilbert}}: A Source Book in the Foundations
  of Mathematics}, vol.~2.
\newblock {Oxford}: {Oxford University Press}.

\bibitem[{Feferman(2003)}]{Feferman2003}
Feferman, Solomon, 2003.
\newblock Introductory note to the {{G\"odel-Bernays}} correspondence.
\newblock In  \parencite{Godel2003}, 41--78.

\bibitem[{Fraenkel(2016)}]{Fraenkel2016}
Fraenkel, Abraham~A., 2016.
\newblock \emph{Recollections of a {{Jewish}} Mathematician in {{Germany}}},
  ed. Jiska {Cohen-Mansfield}.
\newblock {Basel}: {Birkh\"auser}.
\newblock \doi{10.1007/978-3-319-30847-0}.

\bibitem[{Fraenkel(1922)}]{Fraenkel1922a}
Fraenkel, Adolf~Abraham, 1922.
\newblock Zu den {{Grundlagen}} der {{Cantor-Zermeloschen Mengenlehre}}.
\newblock \emph{Mathe\-ma\-tische An\-nalen}~86: 230--237.
\newblock \doi{10.1007/BF01457986}.

\bibitem[{Fraenkel(1926)}]{Fraenkel1926}
Fraenkel, Adolf~Abraham, 1926.
\newblock {Axiomatische Theorie der geordneten Mengen. (Untersuchungen \"uber
  die Grundlagen der Mengenlehre. II.)}.
\newblock \emph{Journal f\"ur die reine und angewandte Mathematik}~155(3):
  129--158.
\newblock \doi{10.1515/crll.1926.155.129}.

\bibitem[{Fraenkel(1928)}]{Fraenkel1928a}
Fraenkel, Adolf~Abraham, 1928.
\newblock {Zusatz zu vorstehendem Aufsatz Herrn v. Neumanns}.
\newblock \emph{Mathematische Annalen}~99(1): 392--393.
\newblock \doi{10.1007/BF01459103}.

\bibitem[{Fraenkel(1937)}]{Fraenkel1937}
Fraenkel, Adolf~Abraham, 1937.
\newblock Ueber eine abgeschwaechte {{Fassung}} des {{Auswahlaxioms}}.
\newblock \emph{The Journal of Symbolic Logic}~2(1): 1--25.
\newblock \doi{10.1017/S002248120003944X}.

\bibitem[{G{\"o}del(1931)}]{Godel1931}
G{\"o}del, Kurt, 1931.
\newblock {\"Uber formal unentscheidbare S\"atze der \emph{Principia
  Mathematica} und verwandter Systeme I}.
\newblock \emph{Monatshefte f\"ur Mathematik und Physik}~38: 173--198.
\newblock \doi{10.1007/BF01700692}.
\newblock Reprinted and translated in \cite[144–195]{Godel1986}.

\bibitem[{G{\"o}del(1986)}]{Godel1986}
G{\"o}del, Kurt, 1986.
\newblock \emph{Publications 1929\textendash 1936}, \emph{Collected Works},
  vol.~1, eds. Solomon Feferman et~al.
\newblock {Oxford}: {Oxford University Press}.

\bibitem[{G{\"o}del(2003{\natexlab{a}})}]{Godel2003}
G{\"o}del, Kurt, 2003{\natexlab{a}}.
\newblock \emph{Correspondence {{A}}\textendash{{G}}}, \emph{Collected Works},
  vol.~4, eds. Solomon Feferman et~al.
\newblock {Oxford}: {Oxford University Press}.

\bibitem[{G{\"o}del(2003{\natexlab{b}})}]{Godel2003a}
G{\"o}del, Kurt, 2003{\natexlab{b}}.
\newblock \emph{Correspondence {{H}}\textendash{{Z}}}, \emph{Collected Works},
  vol.~5, eds. Solomon Feferman et~al.
\newblock {Oxford}: {Oxford University Press}.

\bibitem[{G{\"o}del(2019)}]{Godel2019}
G{\"o}del, Kurt, 2019.
\newblock \emph{{Philosophie I Maximen 0 / Philosophy I Maxims 0}}, ed.
  Eva-Maria Engelen.
\newblock No.~1 in {Philosophical notebooks}. {Berlin}: {De Gruyter}.
\newblock \doi{10.1515/9783110585605}.

\bibitem[{{Grattan-Guinness}(1979)}]{Grattan-Guinness1979}
{Grattan-Guinness}, Ivor, 1979.
\newblock In memoriam {{Kurt G\"odel}}: His 1931 correspondence with
  {{Zermelo}} on his incompletability theorem.
\newblock \emph{Historia Mathematica}~6(3): 294--304.
\newblock \doi{10.1016/0315-0860(79)90127-7}.

\bibitem[{{Grattan-Guinness}(2000)}]{Grattan-Guinness2000}
{Grattan-Guinness}, Ivor, 2000.
\newblock \emph{The Search for Mathematical Roots 1870\textendash 1940: Logics,
  Set Theories and the Foundations of Mathematics from {{Cantor}} through
  {{Russell}} to {{G\"odel}}}.
\newblock {Princeton, N.J.}: {Princeton University Press}.

\bibitem[{Hashagen(2006{\natexlab{a}})}]{Hashagen2006}
Hashagen, Ulf, 2006{\natexlab{a}}.
\newblock {Johann Ludwig Neumann von Margitta (1903\textendash 1957) Teil 1:
  Lehrjahre eines j\"udischen Mathematikers w\"ahrend der Zeit der Weimarer
  Republik}.
\newblock \emph{Informatik-Spektrum}~29(2): 133--141.
\newblock \doi{10.1007/s00287-006-0072-1}.

\bibitem[{Hashagen(2006{\natexlab{b}})}]{Hashagen2006a}
Hashagen, Ulf, 2006{\natexlab{b}}.
\newblock {Johann Ludwig Neumann von Margitta (1903\textendash 1957) Teil 2:
  Ein Privatdozent auf dem Weg von Berlin nach Princeton}.
\newblock \emph{Informatik-Spektrum}~29(3): 227--236.
\newblock \doi{10.1007/s00287-006-0084-x}.

\bibitem[{Hashagen(2010)}]{Hashagen2010}
Hashagen, Ulf, 2010.
\newblock Die {{Habilitation}} von {{John}} von {{Neumann}} an der
  {{Friedrich-Wilhelms-Universit\"at}} in {{Berlin}}: {{Urteile}} \"uber einen
  ungarisch-j\"udischen {{Mathematiker}} in {{Deutschland}} im {{Jahr}} 1927.
\newblock \emph{Historia Mathematica}~37(2): 242--280.
\newblock \doi{10.1016/j.hm.2009.04.002}.

\bibitem[{Herbrand(1930)}]{Herbrand1930}
Herbrand, Jacques, 1930.
\newblock {Recherches sur la th\'eorie de la d\'emonstration}.
\newblock Ph.D. thesis, University of Paris.
\newblock Reprinted in \cite[36--153]{Herbrand1968}. English translation in
  \cite[44--202]{Herbrand1971}.

\bibitem[{Herbrand(1931)}]{Herbrand1931}
Herbrand, Jacques, 1931.
\newblock {Sur la non-contradiction de l'arithm\'etique}.
\newblock \emph{Journal f\"ur die Reine und Angewandte Mathematik}~166: 1--8.
\newblock \doi{10.1515/crll.1932.166.1}.
\newblock Reprinted in \cite[221--232]{Herbrand1968}. English translation in
  \cite[618--628]{vanHeijenoort1967} and in \cite[282--298]{Herbrand1971}.

\bibitem[{Herbrand(1968)}]{Herbrand1968}
Herbrand, Jacques, 1968.
\newblock \emph{{\'Ecrits logiques}}.
\newblock {Paris}: {Presses universitaires de France}.

\bibitem[{Herbrand(1971)}]{Herbrand1971}
Herbrand, Jacques, 1971.
\newblock \emph{Logical Writings}, ed. Warren~D. Goldfarb.
\newblock {Harvard University Press}.

\bibitem[{Hilbert(1923)}]{Hilbert1923}
Hilbert, David, 1923.
\newblock {Die logischen Grundlagen der Mathematik}.
\newblock \emph{Mathe\-ma\-tische An\-nalen}~88(1\textendash 2): 151--165.
\newblock \doi{10.1007/BF01448445}.
\newblock Translated in \cite[1134--1148]{Ewald1996}.

\bibitem[{Jan{\'e}(2001)}]{Jane2001}
Jan{\'e}, Ignacio, 2001.
\newblock Reflections on {{Skolem}}'s relativity of set-theoretical concepts.
\newblock \emph{Philosophia Mathematica}~9(2): 129--153.
\newblock \doi{10.1093/philmat/9.2.129}.

\bibitem[{Kanamori(2009)}]{Kanamori2009}
Kanamori, Akihiro, 2009.
\newblock Bernays and set theory.
\newblock \emph{Bulletin of Symbolic Logic}~15(1): 43--69.
\newblock \doi{10.2178/bsl/1231081769}.

\bibitem[{Kanamori(2012)}]{Kanamori2012}
Kanamori, Akihiro, 2012.
\newblock In praise of replacement.
\newblock \emph{The Bulletin of Symbolic Logic}~18(1): 46--90.
\newblock \doi{10.2178/bsl/1327328439}.

\bibitem[{Kanckos and Lethen(2021)}]{KanckosLethen2021}
Kanckos, Annika and Tim Lethen, 2021.
\newblock The development of {{G\"odel}}'s ontological proof.
\newblock \emph{The Review of Symbolic Logic}~14(4): 1011--1029.
\newblock \doi{10.1017/S1755020319000479}.

\bibitem[{Kennedy(2020)}]{Kennedy2020}
Kennedy, Juliette, 2020.
\newblock \emph{G\"odel, {{Tarski}} and the Lure of Natural Language: Logical
  Entanglement, Formalism Freeness}.
\newblock {Cambridge}: {Cambridge University Press}.
\newblock \doi{10.1017/9780511998393}.

\bibitem[{Lethen(2021)}]{Lethen2021}
Lethen, Tim, 2021.
\newblock \emph{{Gespr\"ache, Vortr\"age, S\'eancen: Kurt G\"odels Wiener
  Protokolle 1937/38: Transkriptionen und Kommentare}}.
\newblock No.~31 in {Ver\"offentlichungen des Instituts Wiener Kreis}. {Cham}:
  {Springer}.
\newblock \doi{10.1007/978-3-030-67606-3}.

\bibitem[{Mancosu(1999)}]{Mancosu1999a}
Mancosu, Paolo, 1999.
\newblock Between {{Vienna}} and {{Berlin}}: The immediate reception of
  {{G\"odel}}'s incompleteness theorems.
\newblock \emph{History and Philosophy of Logic}~20: 33--45.
\newblock \doi{10.1080/014453499298174}.

\bibitem[{Mancosu(2004)}]{Mancosu2004}
Mancosu, Paolo, 2004.
\newblock Review of {{Kurt G\"odel}}, {{{\emph{Collected Works}}}}, vols.
  {{IV}} and {{V}}, {{Solomon Feferman}}, et al., eds. {{Oxford}}: {{Oxford
  University Press}}, 2003.
\newblock \emph{Notre Dame Journal of Formal Logic}~45: 109--125.
\newblock \doi{10.1305/ndjfl/1095386647}.

\bibitem[{Meschkowski(1967)}]{Meschkowski1967}
Meschkowski, Herbert, 1967.
\newblock \emph{{Probleme des Unendlichen. Werk und Leben Georg Cantors}}.
\newblock {Braunschweig}: {Vieweg}.

\bibitem[{Mirimanoff(1917)}]{Mirimanoff1917}
Mirimanoff, Dmitri, 1917.
\newblock Les antinomies de {{Russell}} et de {{Burali-Forti}} et le probl\`eme
  fondamental de la th\'eorie des ensembles.
\newblock \emph{L'Enseignement Math\'ematique}~19(1\textendash 2): 37--52.
\newblock \doi{10.5169/seals-17315}.

\bibitem[{Peckhaus(1990)}]{Peckhaus1990a}
Peckhaus, Volker, 1990.
\newblock `{{Ich}} habe mich wohl geh\"utet, alle {{Patronen}} auf einmal zu
  verschie\ss en'. {{Ernst Zermelo}} in {{G\"ottingen}}.
\newblock \emph{History and Philosophy of Logic}~11(1): 19--58.
\newblock \doi{10.1080/01445349008837156}.

\bibitem[{Peckhaus(1992)}]{Peckhaus1992}
Peckhaus, Volker, 1992.
\newblock Hilbert, {{Zermelo}} und die {{Institutionalisierung}} der
  mathematischen {{Logik}} in {{Deutschland}}.
\newblock \emph{Berichte zur Wissenschaftsgeschichte}~15(1): 27--38.
\newblock \doi{10.1002/bewi.19920150107}.

\bibitem[{Peckhaus(2005)}]{Peckhaus2005}
Peckhaus, Volker, 2005.
\newblock Pro and contra {{Hilbert}}: {{Zermelo}}'s set theories.
\newblock \emph{Philosophia Scienti\ae. Travaux d'histoire et de philosophie
  des sciences}~CS 5: 199--215.
\newblock \doi{10.4000/philosophiascientiae.390}.

\bibitem[{Presburger(1930)}]{Presburger1930}
Presburger, Moj{\.e}sz, 1930.
\newblock {\"Uber die Vollst\"andigkeit eines gewisse Systems der Arithmetik
  ganzer Zahlen, in welchem die Addition als einzige Operation hervortritt}.
\newblock In \emph{{Comptes-rendus du I Congr\`es des Math\'ematiciens des Pays
  Slaves, Varsovie 1929}}, ed. Franciszek Leja, 92--101.
\newblock English translations in \cite{Stansifer1984} and
  \cite{Presburger1991}.

\bibitem[{Presburger(1991)}]{Presburger1991}
Presburger, Moj{\.z}esz, 1991.
\newblock On the completeness of a certain system of arithmetic of whole
  numbers in which addition occurs as the only operation.
\newblock \emph{History and Philosophy of Logic}~12(2): 225--233.
\newblock \doi{10.1080/014453409108837187}.

\bibitem[{Sieg(2003)}]{Sieg2003}
Sieg, Wilfried, 2003.
\newblock Introductory note to the {{G\"odel-Herbrand}} correspondence.
\newblock In  \parencite{Godel2003a}, 3--13.

\bibitem[{Skolem(1923)}]{Skolem1923}
Skolem, Thoralf, 1923.
\newblock {Einige Bemerkungen zur axiomatischen Begr\"undung der Mengenlehre}.
\newblock In \emph{{Matematikerkongressen i Helsingfors 4--7 Juli 1922. Den
  femte skandinaviska matematikerkongressen, Redog\"orelse}}, 217--232.
  {Helsinki}: {Akademiska Bokhandeln}.
\newblock Reprinted in \cite[137--152]{Skolem1970}; English translation in
  \cite[290--301]{vanHeijenoort1967}.

\bibitem[{Skolem(1970)}]{Skolem1970}
Skolem, Thoralf, 1970.
\newblock \emph{Selected Works in Logic}, ed. Jens~Erik Fenstad.
\newblock {Oslo}: {Universitetsforlaget}.

\bibitem[{Stansifer(1984)}]{Stansifer1984}
Stansifer, Ryan, 1984.
\newblock Presburger's article on integer airthmetic: Remarks and translation.
\newblock Tech. Rep. TR84-639, {Cornell University, Computer Science
  Department}.
\newblock \urlprefix\url{https://hdl.handle.net/1813/6478}.

\bibitem[{{van Atten}(2015)}]{vanAtten2015}
{van Atten}, Mark, 2015.
\newblock \emph{Essays on {{G\"odel}}'s Reception of {{Leibniz}}, {{Husserl}},
  and {{Brouwer}}}.
\newblock No.~35 in Logic, {{Epistemology}}, and the {{Unity}} of {{Science}}.
  {Cham}: {Springer}.
\newblock \doi{10.1007/978-3-319-10031-9}.

\bibitem[{{van Heijenoort}(1967)}]{vanHeijenoort1967}
{van Heijenoort}, Jean, ed., 1967.
\newblock \emph{From {{Frege}} to {{G\"odel}}: A Source Book in Mathematical
  Logic, 1897\textendash 1931}.
\newblock {Cambridge, MA}: {Harvard University Press}.

\bibitem[{{von Neumann}(1925)}]{vonNeumann1925}
{von Neumann}, Johann, 1925.
\newblock {Eine Axiomatisierung der Mengenlehre}.
\newblock \emph{Journal f\"ur die reine und angewandte Mathematik}~154(4):
  219--240.
\newblock \doi{10.1515/crll.1925.154.219}.
\newblock Reprinted in \cite[34--56]{vonNeumann1961}. English translation in
  \cite[393--413]{vanHeijenoort1967}.

\bibitem[{{von Neumann}(1927)}]{vonNeumann1927}
{von Neumann}, Johann, 1927.
\newblock {Zur Hilbertschen Beweistheorie}.
\newblock \emph{Mathematische Zeitschrift}~26(1): 1--46.
\newblock \doi{10.1007/BF01475439}.
\newblock Reprinted in \cite[256--300]{vonNeumann1961}.

\bibitem[{{von Neumann}(1928{\natexlab{a}})}]{vonNeumann1928}
{von Neumann}, Johann, 1928{\natexlab{a}}.
\newblock {Die Axiomatisierung der Mengenlehre}.
\newblock \emph{Mathe\-mati\-sche Zeit\-schrift}~27: 669--752.
\newblock \doi{10.1007/BF01171122}.
\newblock Reprinted in \cite[339--422]{vonNeumann1961}.

\bibitem[{{von Neumann}(1928{\natexlab{b}})}]{vonNeumann1928a}
{von Neumann}, Johann, 1928{\natexlab{b}}.
\newblock {\"Uber die Definition durch transfinite Induktion und verwandte
  Fragen der allgemeinen Mengenlehre}.
\newblock \emph{Mathematische Annalen}~99(1): 373--391.
\newblock \doi{10.1007/BF01459102}.
\newblock Reprinted in \cite[320--338]{vonNeumann1961}.

\bibitem[{{von Neumann}(1929)}]{vonNeumann1929a}
{von Neumann}, Johann, 1929.
\newblock {\"Uber eine Widerspruchfreiheitsfrage in der axiomatischen
  Mengenlehre.}
\newblock \emph{Journal f\"ur die reine und angewandte Mathematik}~160(4):
  227--241.
\newblock \doi{10.1515/crll.1929.160.227}.
\newblock Reprinted in \cite[494--508]{vonNeumann1961}.

\bibitem[{{von Neumann}(1961)}]{vonNeumann1961}
{von Neumann}, John, 1961.
\newblock \emph{Logic, Theory of Sets and Quantum Mechanics}, ed.
  Abraham~Haskel Taub.
\newblock No.~1 in Collected Works. {Oxford}: {Pergamon Press}.

\bibitem[{{von Plato}(2020)}]{vonPlato2020}
{von Plato}, Jan, ed., 2020.
\newblock \emph{Can Mathematics Be Proved Consistent? {{G\"odel}}'s Shorthand
  Notes \& Lectures on Incompleteness}.
\newblock Sources and {{Studies}} in the {{History}} of {{Mathematics}} and
  {{Physical Sciences}}. {Cham}: {Springer}.
\newblock \doi{10.1007/978-3-030-50876-0_1}.

\bibitem[{Weyl(1918)}]{Weyl1918}
Weyl, Hermann, 1918.
\newblock \emph{{Das Kontinuum}}.
\newblock {Leipzig}: {Veit}.

\bibitem[{Zach(2003)}]{Zach2003}
Zach, Richard, 2003.
\newblock The practice of finitism: Epsilon calculus and consistency proofs in
  {{Hilbert}}'s {{Program}}.
\newblock \emph{Synthese}~137(1/2): 211--259.
\newblock \doi{10.1023/A:1026247421383}.

\bibitem[{Zach(2004)}]{Zach2004a}
Zach, Richard, 2004.
\newblock Hilbert's `{{{\emph{Verungl\"uckter Beweis}}}}', the first epsilon
  theorem, and consistency proofs.
\newblock \emph{History and Philosophy of Logic}~25(2): 79--94.
\newblock \doi{10.1080/01445340310001606930}.

\bibitem[{Zermelo(1932)}]{Zermelo1932}
Zermelo, Ernst, 1932.
\newblock {\"Uber Stufen der Quantifikation und die Logik des Unendlichen}.
\newblock \emph{Jahresbericht der Deutschen Mathematiker-Vereinigung}~41:
  85--88.

\end{thebibliography}

\end{document}